\documentclass[final]{siamltex}

\usepackage{amsmath}
\usepackage{amssymb}
\usepackage{graphicx}
\usepackage{caption}
\usepackage{subcaption}
\usepackage{color}

\newcommand{\half}{\frac{1}{2}} % 1/2
\newcommand{\area}{\mathcal{A}}

\renewcommand{\vec}[1]{\mathbf{#1}} % Vector (bold)
\newcommand{\labelEq}[1]{\label{eq:#1}} % Use this to label an equation.
\newcommand{\refEq}[1]{(\ref{eq:#1})}   % Use this to reference an equation.
\newcommand{\labelSec}[1]{\label{sec:#1}} % Use this to label a section.
\newcommand{\refSec}[1]{\S\ref{sec:#1}} % Use this to reference a section.
\newcommand{\ib}{{\boldsymbol{i}}}    % bold italic i
\newcommand{\ub}{{\boldsymbol{u}}}    % bold italic u
\newcommand{\xb}{{\boldsymbol{x}}}    % bold italic x
\newcommand{\pb}{{\boldsymbol{p}}}    % (bold italic) multi-index p
\newcommand{\avg}[1]{\left<#1\right>} % Cell/face averaged values
\newcommand{\tvol}{\text{vol}}

\newcommand{\nbold}{{\boldsymbol{n}}}

\newcommand{\ed}{{\boldsymbol{e}^d}}

\newcommand{\xbold}{{\boldsymbol{x}}}
\newcommand{\ybold}{{\boldsymbol{y}}}
\newcommand{\pbold}{{\boldsymbol{p}}}

\newcommand{\vol}{\mathcal{V}}

\newcommand{\zerobold}{{\boldsymbol{0}}}
\newcommand{\xz}{{\boldsymbol{x}_0}}

\title{A higher-order finite-volume discretization method \\
       for Poisson's equation in cut cell geometries 
       }

\author{ D. Devendran  \footnotemark[2], 
         D. T. Graves  \footnotemark[2], 
         H. Johansen   \footnotemark[2] }

\footnotetext[2]{Lawrence Berkeley National Laboratory, Berkeley, CA. Research at LBNL was supported financially by the Office of Advanced Scientific Computing
Research of the US Department of Energy under contract number
DE-AC02-05CH11231.}

\begin{document}

\maketitle

\begin{abstract}
  We present a method for generating higher-order finite volume
  discretizations for Poisson's equation on Cartesian cut cell grids
  in two and three dimensions.  The discretization is in
  flux-divergence form, and stencils for the flux are computed by
  solving small weighted least-squares linear systems.  Weights are
  the key in generating a stable discretization.  We apply the method to solve
  Poisson's equation on a variety of geometries, and we demonstrate
  that the method can achieve second and fourth order accuracy in both
  truncation and solution error for these examples.  We also show that
  the Laplacian operator has only stable eigenvalues for each of these
  examples.
\end{abstract}


\pagestyle{myheadings}
\thispagestyle{plain}
\markboth{D. DEVENDRAN, D.T. GRAVES, H. JOHANSEN}{HIGHER-ORDER CUT
  CELL METHOD FOR POISSON}

\newpage
\section{Introduction} 

There are many numerical approaches to solve Poisson's equation in
complex geometries.  Green function approaches \cite{McKenneyFMM1994,
  Greengard1996, Cheng1999}, such as the fast multipole method, are
fast and near-optimal in complexity, but they are not conservative.
Also, they cannot be easily extended to variable and tensor
coefficient Poisson operators, which are important in the earth
sciences and multi-material problems.

  Another popular approach is to use the
finite element method, which has a number of advantages.  These
advantages include negative-definite discrete operators, higher-order
accuracy, and ease of extension to variable coefficients.  The
conditioning and accuracy of the discrete finite element operator can
be strongly mesh-dependent, however \cite{BRENNERFEMBOOK}.
Unfortunately, generating meshes with higher-order conforming elements
for complex 3D domains is still an expensive, globally-coupled
computation, and an open area of research \cite{PirzadehFEMGridGen2010}.

This motivates the need for simpler grid generation.  Cut cells are a
simple way of addressing this.  In a cut cell (or embedded boundary)
method, the discrete domain is the intersection of the complex
geometry with a regular Cartesian grid.  Such intersections are local,
and can be calculated very efficiently in parallel, enabling fast
computation of solution-dependent moving boundaries
\cite{AftosmisBergerMelton, EBGeometryPaper}.  The complexity of
dealing with complex geometries is shifted back to the discretization
approach. The cut-cell approach has been used successfully to solve
Poisson's equation in finite volume \cite{johansenColella:1998,
  schwartzETAL:2006} and finite difference \cite{GibouFedkiw2005,
  LevequeLing1994} discretizations.

For many problems, such as heat and mass transfer, discrete
conservation is important.  Finite volume methods are discretely
conservative by construction because they are in discrete
flux-divergence form \cite{LEVEQUEBOOK}.  Previous finite volume
methods for Poisson's equation are first order in truncation error
near the embedded boundary and second order in solution error
\cite{johansenColella:1998, schwartzETAL:2006}.

We present a method for generating higher-order finite volume
discretizations for Poisson's equation on Cartesian cut cell grids in
two and three dimensions.  The discretization is in flux-divergence
form.  We compute stencils for the flux by solving small weighted
least squares systems.  In principle, the method can produce
discretizations for any given order of accuracy.  In Section
\refSec{results}, we apply the method to solve Poisson's equations on
a number of geometries, and we demonstrate the method can achieve both
second order and fourth order convergence in both truncation and
solution error.

This paper is organized as follows.  In Section \refSec{method}, we
introduce our method.  In Section \refSec{results}, we present 2D and
3D examples that demonstrate convergence with grid refinement.  We
also show that the Laplacian operator has only stable eigenvalues.
Finally, we demonstrate that the method is robust under small
perturbations in the geometry.

\section{Method}
\labelSec{method}

%{\tt method.tex} \\

We design a conservative finite volume method to solve Poisson's equation
\begin{equation*}
 \Delta \phi = \rho
\end{equation*}
for the potential $\phi$ on a domain $\Omega$ with a charge distribution $\rho$.  First, we write the equation in flux-divergence form:
\begin{equation*}
	\nabla \cdot \nabla \phi = \rho .
\end{equation*}
Integrating over an arbitrary region $V \subseteq \Omega$ and applying the divergence theorem gives
\begin{equation}
\labelEq{laplacian over V}
	\int_{V} \nabla \cdot \nabla \phi d\vol= \int_{\partial V} \nabla \phi \cdot \nbold d\area.
\end{equation}
Our method is based on using a higher-order interpolant of $\phi$ to approximate the flux $\int_{\partial V} \nabla \phi \cdot \nbold d\area$.

\subsection{Spatial Notation}

\quad 
Our computational domain is a set of distinct, contiguous volumes,
  $\{\vol_v\}$, each of which is part of an intersection of $\Omega$ 
  with a cell $\vol_\ib$ in a regular grid of grid spacing $h$,
\begin{equation}
\nonumber
\vol_\ib = [i_1 h, (i_1 + 1)h] \otimes ... \otimes [i_D h, (i_D + 1)h] 
      \equiv [\ib h, (\ib + \ub)h] ,
\end{equation}
  where the index $\ib = (i_1, ..., i_D)$, and $\ub = (1, ..., 1)$.
Note that we use the index $v$ to uniquely identify a volume;
  for a given regular cell $\vol_\ib$, there may be more than one 
  $\vol_v$ such that 
  $\cup \{\vol_v \cap \vol_\ib \} = \vol_\ib \cap \Omega$,
  especially in the case of very complex geometries.

The grid-aligned faces associated with $\vol_\ib$ in the $\pm d$ directions are 
  identified by an additional half index, $\ib \pm \half\vec{e}^d$, where
  $\vec{e}^d$ is the unit vector with components $e^d_i = 1$ if $i = d,\, 0$
  otherwise.
For example, 
\begin{equation}
\nonumber
  \area_{\ib + \half\vec{e}^d} 
  \equiv \left[\ib h, (\ib + \ub - \vec{e}^d)h \right] \, .
\end{equation}
For a given volume $\vol_v$, its surface $\partial \vol_v$
  is discretized into grid-aligned faces, 
  $\area_{v \pm \half\vec{e}^d} = \area_{\ib \pm \half\vec{e}^d} \cap \partial \vol_v$,
  which are shared between neighboring volumes.
$\vol_v$ may also contain a portion of the domain boundary, which
  we indicate with 
  $\area^b_v = \partial \Omega \cap \partial \vol_v$.
In either case, we use the index $f$ to provide a unique global index into the
  set of all such faces, $\{f(v): \area_f \in \partial \vol_v \}$.
When discussing volume- or face-average quantities in the sections below, 
  we will use notational shortcuts where $v$ or $f$ are associated 
  with $\vol_v$ or $\area_f$, respectively.
For example,
\begin{align}
\labelEq{volavg}
  \avg{\phi}_v 
  & \equiv \frac{1}{|\vol_v|} \int_{\vol_v} \phi \, d\vol \hbox{ , and }
\\
%\labelEq{faceavg}
\nonumber
  \avg{\phi}_f 
  & \equiv \frac{1}{|\area_f|} \int_{\area_f} \phi \, d\area \, .
\end{align}

Volumes and faces contained within $\Omega$ that do not contain a portion of the domain boundary are called ``full,'' whereas those that do are ``cut'' by the embedded boundary.  We will often identify irregular faces and volumes in terms of their ``fraction'' of a regular one, that is:
\begin{align}
%\label{eq:volfraction}
\nonumber
    |\vol_v| 
    & = \kappa_v h^D \hbox{ , and }
\\ \nonumber
    |\area_f| 
    & = \alpha_f h^{D-1} \, ,
\end{align}
where $\kappa_v$ and $\alpha_f$ are called the volume- and area-fraction, respectively.

Finally, we define the volume moments and face moments that show up in
our discretization.  In the paper, we use multi-index notation. In
particular, for $\pbold = (p_1, ..., p_D)$,
$\xbold = (x_1, ..., x_D)$, and $\ybold = (y_1, ..., y_D)$, we define
\begin{equation}
\nonumber
 (\xbold + \ybold)^\pbold = \prod_{d=1}^D (x_d + y_d)^{p_d} .
\end{equation} 

The $\pbold$th volume moment of $\vol_v$ is
  \begin{equation}
\labelEq{volume moment}
    m^{\pbold}_v(\xz) = \int_{\vol_v} (\xb - \xz)^{\pbold} d\vol .
  \end{equation}
The $\pbold$th face moment of $\area_{f}$ is
  \begin{equation}
\labelEq{face moment}
   m^{\pbold}_{f}(\xz) = \int_{\area_{f}} (\xb - \xz)^{\pbold} d\area .
  \end{equation}
For an embedded boundary face, it is useful to define a second face moment that includes the normal to the face:
  \begin{equation}
\labelEq{face moment with normal}
   m^{\pbold}_{d, f}(\xz) = \int_{\area_{f}} (\xb - \xz)^{\pbold} n_d(\xb) d\area ,
  \end{equation}
where $n_d$ is the $d$th component of the outward unit normal to $f$.

%%%%%%%
\subsection{Flux Approximation}

Setting $V = \vol_v$ in \refEq{laplacian over V} and dividing by the volume of $\vol_v$ we get
\begin{equation}
\nonumber
 \frac{1}{|\vol_v|} \int_{\vol_v} \nabla \cdot \nabla \phi  \, d\vol = \frac{1}{|\vol_v|} \sum_{f(v)} \int_{\area_f} \, \nabla \phi \cdot \nbold \, d\area \, .
\end{equation}
We approximate the flux by replacing $\phi$ by a polynomial interpolant.  Suppose 
\begin{equation}
\nonumber
\psi(\xb) = \sum_{|\pbold| < P} c_{\pbold} (\xb - \xz)^{\pbold}
\end{equation}
is a polynomial interpolant of $\phi$ such that $\phi(\xb) = \psi(\xz) + O(|\xb - \xz|^P)$.  Then
\begin{equation}
\nonumber
  \int_{\area_f} \nabla \phi \cdot \nbold d\area = \sum_{|\pbold| < P} c_{\pbold} \int_{\area_f} \nabla (\xb - \xz)^{\pbold} \cdot \nbold d\area + O(h^{P+D-1}). 
\end{equation}

\subsubsection{Calculating the Polynomial Interpolant} \label{sec: polynomial interpolant}
To solve for the coefficients $c_\pbold$, we create a system of equations using voume-averaged values of neighboring volumes and face-averaged values of neighboring boundary faces.  For a face $f$, we require 
  \begin{equation}
\nonumber
    \avg{\psi}_v = \avg{\phi}_v ,
  \end{equation}
for all neighboring volumes $\vol_v$.  Using equations \refEq{volavg} and \refEq{volume moment}, this simplifies to
  \begin{equation}
  \labelEq{volume moment eq for cp}
   \frac{1}{|\vol_v|} \sum_{|\pbold| < P} c_{\pbold} m^{\pbold}_v = \avg{\phi}_v .
  \end{equation}

If we are given the Dirichlet boundary condition $\phi = g$ on a neighboring boundary face $f_b$, then we require 
  \begin{equation*}
   \frac{1}{|\area_{f_b}|} \sum_{|\pbold| < P} c_{\pbold} \int_{\area_{f_b}} (\xb - \xz)^{\pbold} d\area = \frac{1}{|\area_{f_b}|} \int_{\area_{f_b}} g d\area ,
  \end{equation*}
or
  \begin{equation}
  \labelEq{Dir moment eq for cp}
   \frac{1}{|\area_{f_b}|} \sum_{|\pbold| < P} c_{\pbold} m^{\pbold}_{f_b}(\xz) = \frac{1}{|\area_{f_b}|} \int_{\area_{f_b}} g d\area ,
  \end{equation}
using equation \refEq{face moment}.

If, instead, the Neumann boundary condition $\nabla \phi \cdot \nbold = g$ is specified on $f_b$ then we require
  \begin{equation*}
   \frac{1}{|\area_{f_b}|} \sum_{d=1}^D \sum_{|\pbold| < P} c_{\pbold} \int_{\area_{f_b}} \partial^{\ed}(\xb - \xz)^{\pbold} n_d d\area = \frac{1}{|\area_{f_b}|} \int_{\area_{f_b}} g d\area ,
  \end{equation*}
or
  \begin{equation}
  \labelEq{Neu moment eq for cp}
   \frac{1}{|\area_{f_b}|} \sum_{d=1}^{D}\sum_{|\pbold| < P} p_d c_{\pbold} m^{\pbold-\ed}_{d, f_b} = \frac{1}{|\area_{f_b}|} \int_{\area_{f_b}} g d\area ,
  \end{equation}
using equation \refEq{face moment with normal}.  Here, we set $m^{\pbold-\ed}_{d, f_b} = 0$ if $\pbold \leq \ed$ to simplify the notation.

Equations \refEq{volume moment eq for cp}, \refEq{Dir moment eq for cp}, and \refEq{Neu moment eq for cp} represent a linear system of equations for the polynomial coefficients $c_{\pbold}$.  Let $c = [\dots \, c_{\pbold} \, \dots]^T$ and
  \begin{equation}
\nonumber
   A_{v \pb} = \frac{m^{\pbold}_v(\xz)}{m^{\zerobold}_v}.
  \end{equation}
Let
  \begin{equation}
\nonumber
   b_{f_b, \pbold} = \frac{m^\pbold_{f_b}}{m^{\zerobold}_{f_b}},
  \end{equation}
if Dirichlet boundary conditions are prescribed on $f_b$ and
  \begin{equation}
\nonumber
   b_{f_b, \pbold} = \frac{1}{m^{\zerobold}_{f_b}} \sum \limits_{d=1}^D p_d \, m^{\pbold - \ed}_{d,f_b}  .
  \end{equation}
if Neumann boundary conditions are prescribed on $f_b$.  Recall, $m^{\zerobold}_v = |\vol_v|$ and $m^{\zerobold}_f = |\area_f|$ (see equations \refEq{volume moment} and \refEq{face moment}).

Finally, let $\Phi_{\tvol} = [\dots \, \avg{\phi}_v \, \dots]^T$ and $G = [\dots \, \avg{g}_{f_b} \, \dots]^T$.  Then the combined linear system is
  \begin{equation}
\nonumber
     \begin{bmatrix} 
    \dots \\
    \dots \, A_{v \pb} \, \dots \\ 
     \dots \\
    \dots \, b_{f_b, \pb} \, \dots \\ 
    \dots
  \end{bmatrix}
    c = 
  \begin{bmatrix}
   \Phi_{\tvol} \\ G 
  \end{bmatrix}
	,
  \end{equation}
or
 \begin{equation}
\nonumber
  Ac = \Phi,
 \end{equation}
where
 \begin{equation}
\labelEq{def of A}
  A = 
     \begin{bmatrix} 
    \dots \\
    \dots \, A_{v \pb} \, \dots \\ 
     \dots \\
    \dots \, b_{f_b, \pb} \, \dots \\ 
    \dots
  \end{bmatrix}
 \end{equation}
and
\begin{equation}
\labelEq{def of phi vec}
\Phi = 
 \begin{bmatrix}
   \Phi_{\tvol} \\ G 
  \end{bmatrix}
.
\end{equation}
If we have sufficiently many neighbors, this is an overdetermined full-rank system.  We use weighted least squares to obtain the coefficients $c$:
\begin{equation}
\nonumber
  c = (WA)^\dagger W\Phi .
\end{equation}
In this paper, we only consider an invertible diagonal weighting matrix $W$.  

\begin{figure}
\begin{subfigure}{0.5\textwidth}
\includegraphics[width=\textwidth]{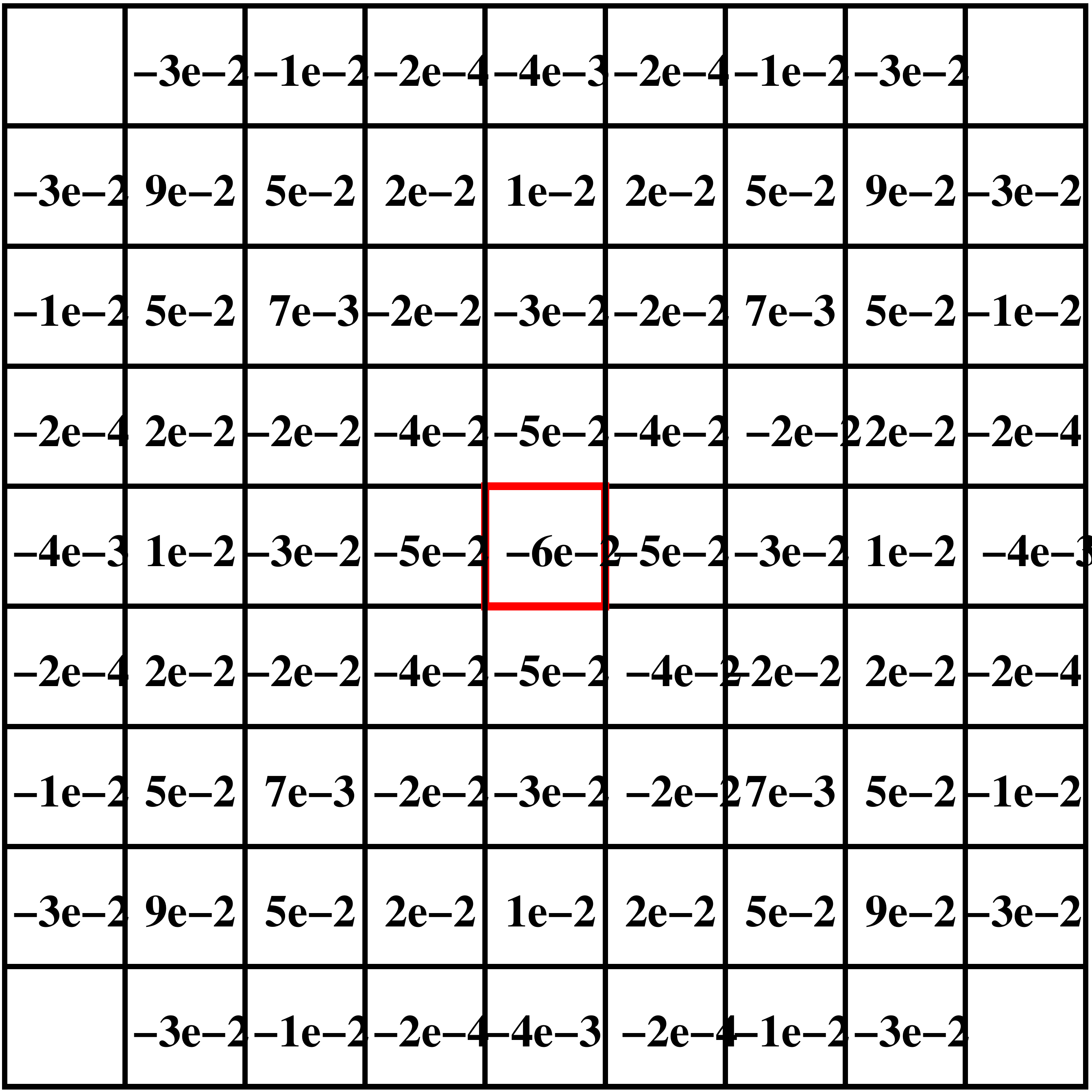}
\caption{Stencil for a regular cell using standard least squares.}
\end{subfigure}
\hspace{0.03\textwidth}
\begin{subfigure}{0.5\textwidth}
\includegraphics[width=\textwidth]{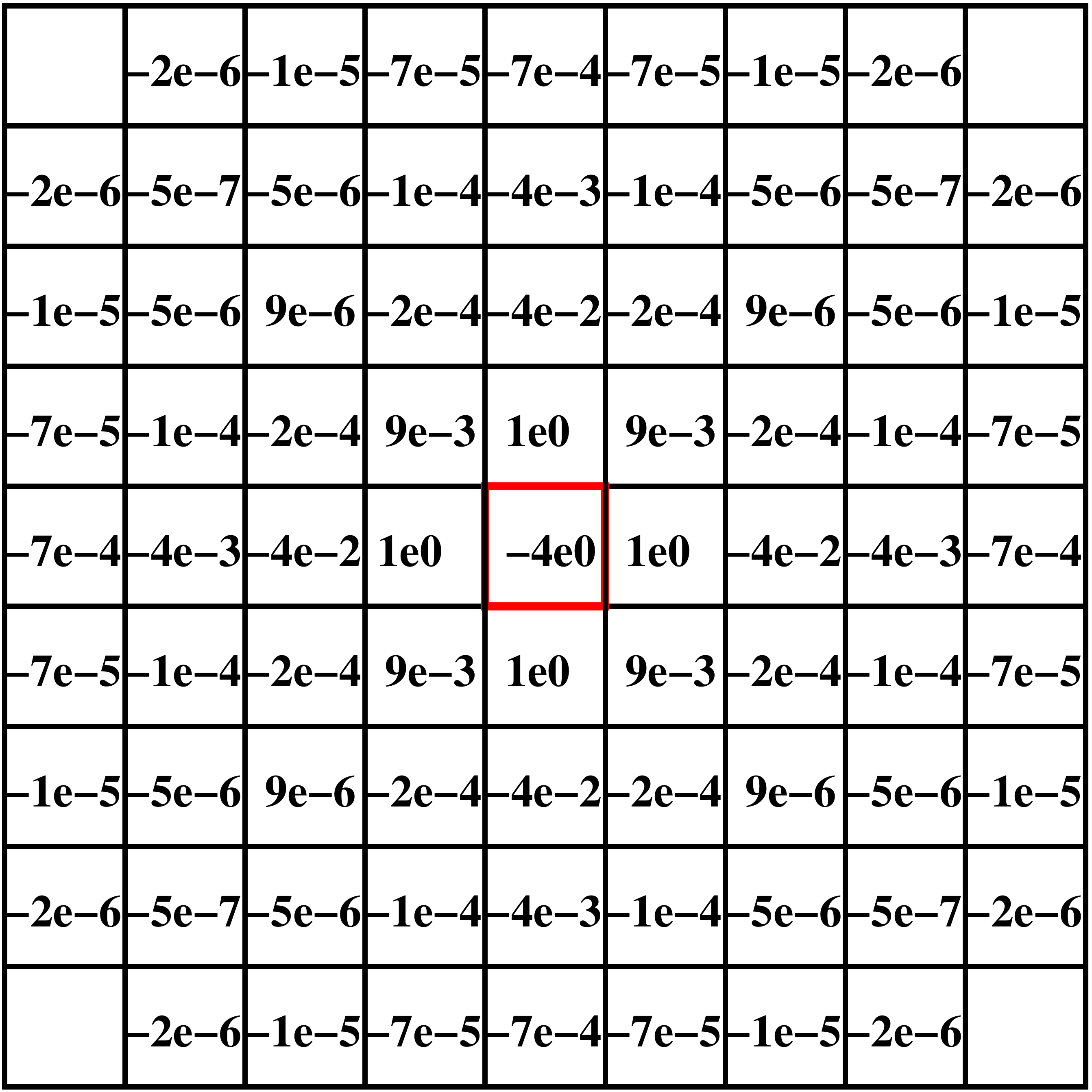}
\caption{Stencil for a regular cell using weighted least squares.}
\end{subfigure}
\begin{subfigure}{0.5\textwidth}
\includegraphics[width=\textwidth]{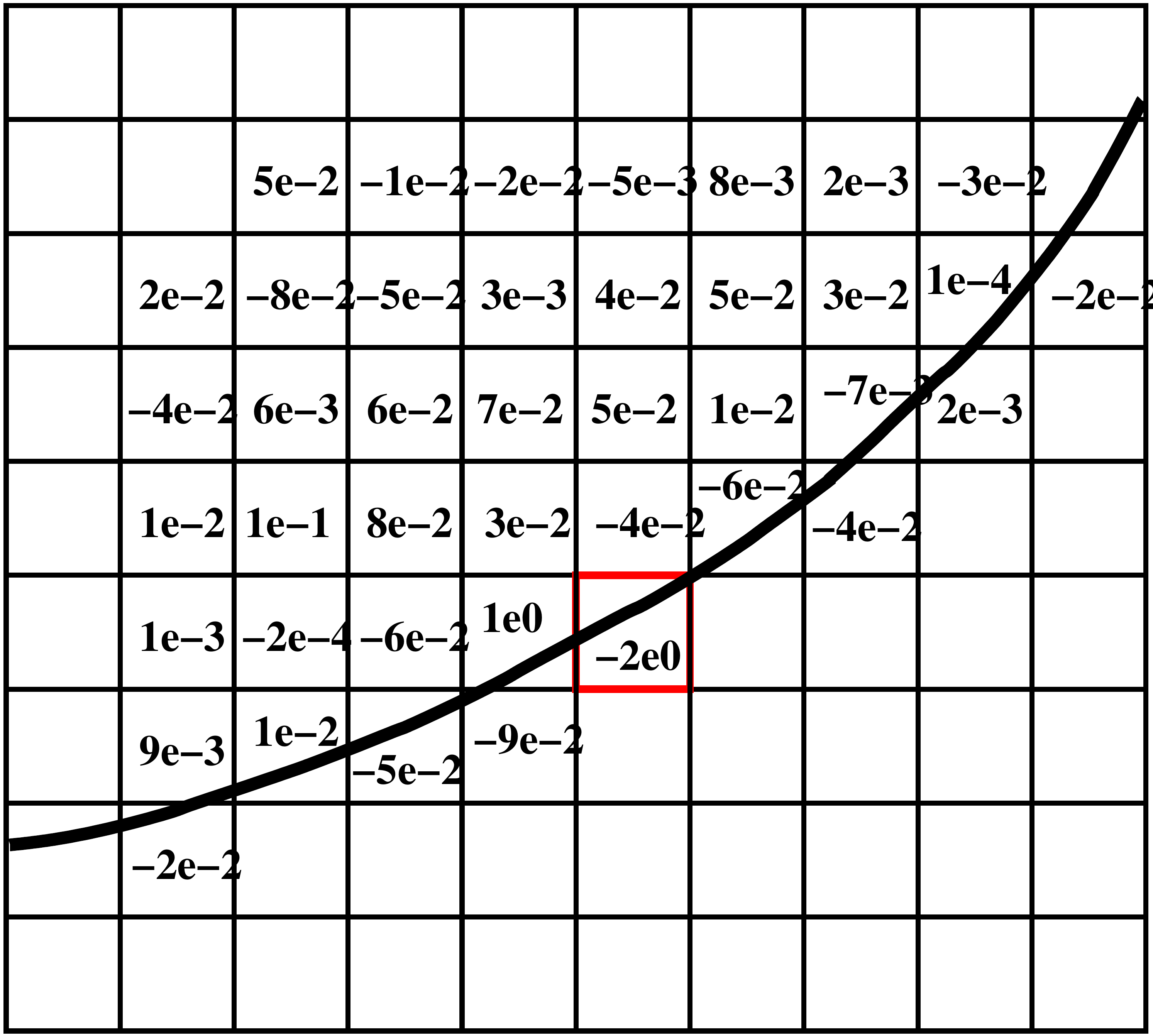}
\caption{Stencil for a cut cell using standard least squares.}
\end{subfigure}
\hspace{0.03\textwidth}
\begin{subfigure}{0.5\textwidth}
\includegraphics[width=\textwidth]{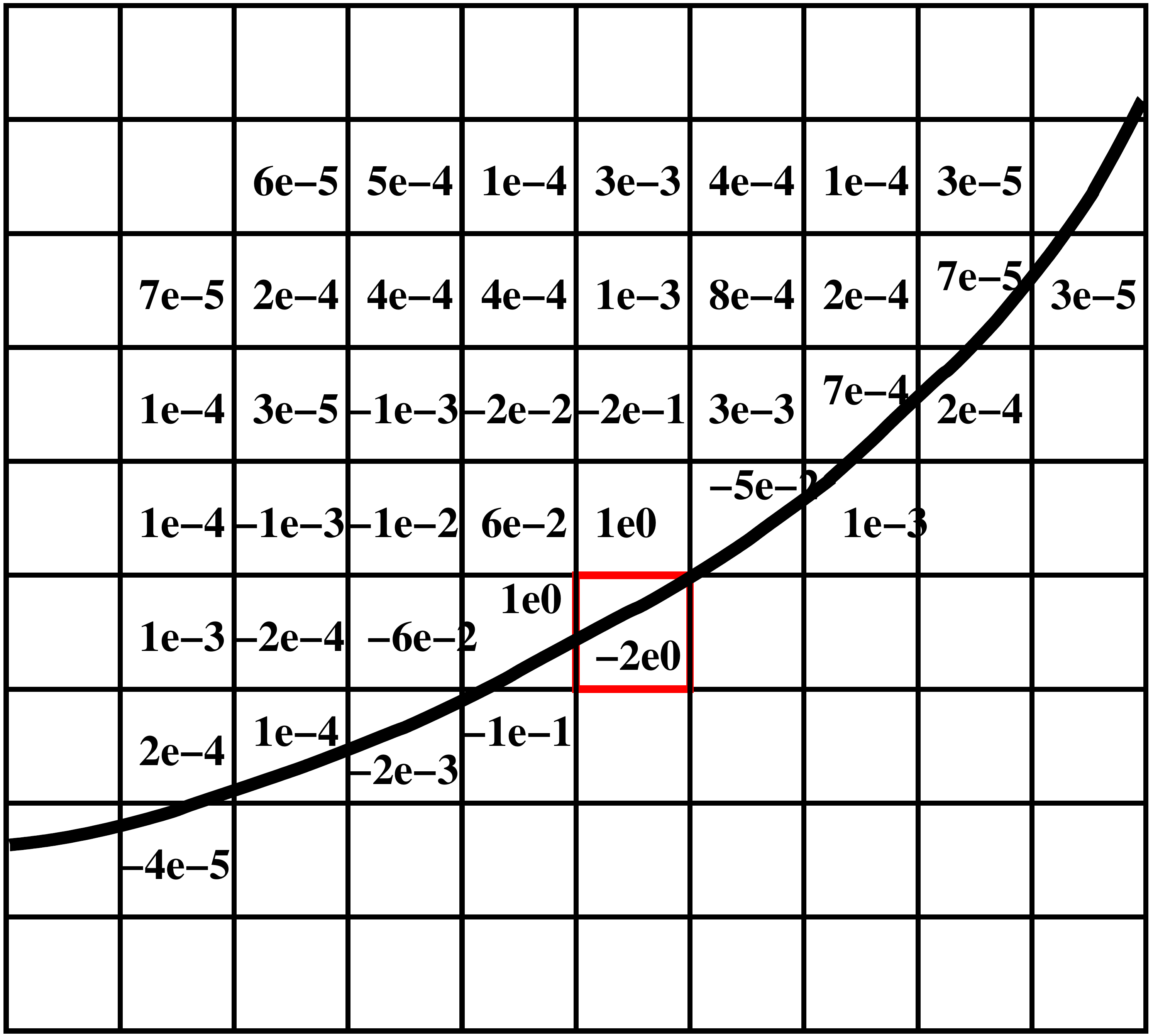}
\caption{Stencil for a cut cell using weighted least squares.}
\end{subfigure}
\caption{Fourth-order Laplacian stencil for the cell highlighted in red.  Without weights, the stencil values do not decay with distance.}
\label{fig:laplacian stencils}
\end{figure}

The weights are extra degrees of freedom in the system, which we use to generate a stable Laplacian operator with eigenvalues that lie on the left half of the complex plane.  Let $F =[\dots \, F_\pbold \, \dots]^T$, where
\begin{equation}
\nonumber
F_{\pbold} =  \int_{\area_f} \nabla (\xb - \xz)^{\pbold} \cdot \nbold d\area .
\end{equation}
Then 
\begin{align}
\nonumber
  \sum_{|\pbold| < P} c_{\pbold} \int_{\area_{f}} (\xb - \xz)^{\pbold} d\area &= F^Tc \\
	 &= F^T(WA)^\dagger W\Phi \nonumber.
\end{align}
The vector 
\begin{equation}
\labelEq{stencil equation solution}
s = W (A^T W)^\dagger F
\end{equation}
is a stencil for computing the flux through $f$.  Note that $s$ satisfies the equation
\begin{equation}
\labelEq{stencil equation}
  A^T s = F ,
\end{equation}
which is an underdetermined system for $s$.  The case
\begin{equation}
\nonumber
 s = (A^T)^\dagger F
\end{equation}
is the solution with minimum $L_2$ norm.  This flux stencil does not
decay with distance from the face (see Figure \ref{fig:laplacian
  stencils}), and it produces an unstable Laplacian operator with
large positive eigenvalues (data not shown).  Equation \refEq{stencil equation solution} is the solution that minimizes $\| W^{-1} s\|_2$.  By choosing weights that decay with distance, we can force the flux stencil to also decay with distance (see Figure \ref{fig:laplacian stencils}).  Our particular choice of weights is given in Section \ref{sec: neighs and weights}.

\subsection{Order of Accuracy}
In general, we cannot calculate the moments exactly.  The geometry and
the corresponding moments are constructed using the method in
\cite{EBGeometryPaper}.  Let $M_{d,f}^{\pbold}(\xz)$ denote the
approximation to the exact moment $m_{d,f}^{\pbold}(\xz)$, so that
\begin{equation}
\nonumber
M_{d,f}^{\pbold}(\xz) = m_{d,f}^{\pbold}(\xz) + O(h^R) .
\end{equation}

Recall that
\begin{equation}
\nonumber
\phi(\xb) = \sum_{|\pbold| < P } c_{\pbold} (\xb - \xz)^{\pbold} + O(|\xb - \xz|^P) .
\end{equation}

In this subsection, we show that the truncation error for the kappa-weighted Laplacian operator is $O(h^{P-2})$ if $R = P+D-2$.

We assume that for face $f$, $\xz$ is some point such that $|\xb - \xz| \leq Ch$ for all points $\xb$ on the face.  Then
\begin{align*}
	\int_{\area_f} \nabla \phi \cdot \nbold d\area &= \sum_{|\pbold| < P} c_{\pbold} \int_{\area_f} \nabla (\xb - \xz)^{\pbold} \cdot \nbold d\area + O(h^{P+D - 2}) \\ 
&= \sum_{d=1}^D \sum_{|\pbold| < P}c_{\pbold} \int_{\area_{f}} p_d (\xb - \xz)^{\pbold-\ed} n_d(\xb) d\area + O(h^{P+D - 2}) \\
&= \sum_{d=1}^D \sum_{|\pbold| < P}p_d c_{\pbold} m_{d,f}^{\pbold-\ed}(\xz) + O(h^{P+D - 2}) \\
	&= \sum_{d=1}^D\sum_{|\pbold| < P}p_d c_{\pbold} M_{d,f}^{\pbold-\ed}(\xz) + O(h^R) + O(h^{P+D - 2}) . \\
\end{align*}

We get
\begin{align*}
\frac{\kappa_v}{|\vol_v|} \int_{\vol_v} \nabla \cdot \nabla \phi  \, d\vol &= \frac{1}{h^D}\sum_{f(v)} \int_{\area_f} \nabla \phi \cdot \nbold d\area \\
			                   &= \frac{1}{h^D} \sum_{f(v)} \sum_{d=1}^D\sum_{|\pbold| < P}p_d c_{\pbold} M_{d,f}^{\pbold-\ed}(\xz) + O(h^R) + O(h^{P+D - 2}) . \\
\end{align*}

At volumes $\vol_v$ that are sufficiently far from the boundary (i.e. they only include regular volumes in their Laplacian stencils), the truncation error is $O(h^{P-1})$ if $R = P+D-1$.  The $O(h^{P-2})$ term in the polynomial interpolation error cancels out because of symmetry in the flux stencils.

The truncation error is not the same for different choices of norms
because of this disparity in the order of accuracy between volumes
near the boundary and interior volumes.  The error in the $L_\infty$ norm is $O(h^{P-2})$.  In the $L_1$ norm, however, it is $O(h^{P-1})$.  If $N = 1/h$, then the number of volumes near the boundary is $O(N)$ whereas the number of interior volumes is $O(N^2)$.  So, the $L_1$ error looks like
	\begin{equation}
 \nonumber
		\sum |e_{ij}|\Delta V \approx N^2(h^{P-1})(h^2) + N(h^{P-2})(h^2) = h^{P-1} .
	\end{equation}

In contrast, the solution error has the same order of accuracy in all
the norms.  Using potential theory arguments (see
\cite{Johansen:1998:EBPoisson} for details), it can be shown that the error is $O(h^{P})$ near the boundary for Dirichlet boundary conditions and $O(h^{P-1})$ for Neumann boundary conditions.  In the interior, the solution error is $O(h^{P-1})$.  So, the order of accuracy in the solution is $Q$ if $P = Q+1$ and $R = Q+D$.   

\subsection{Neighbors and Weights} \label{sec: neighs and weights}

Recall, the polynomial interpolant $\psi$ is 
\begin{equation}
 \nonumber
	\psi = \sum_{|\pbold| < P} c_{\pbold} (\xb - \xz)^{\pbold}.
\end{equation}
In our method, we choose $\xz$ to be the face center when $f$ is a grid-aligned face. If $f$ is an embedded boundary face, $\xz$ is the center of the cell cut by $f$.

\begin{figure}
\begin{subfigure}{0.5\textwidth}
\includegraphics[width=\textwidth]{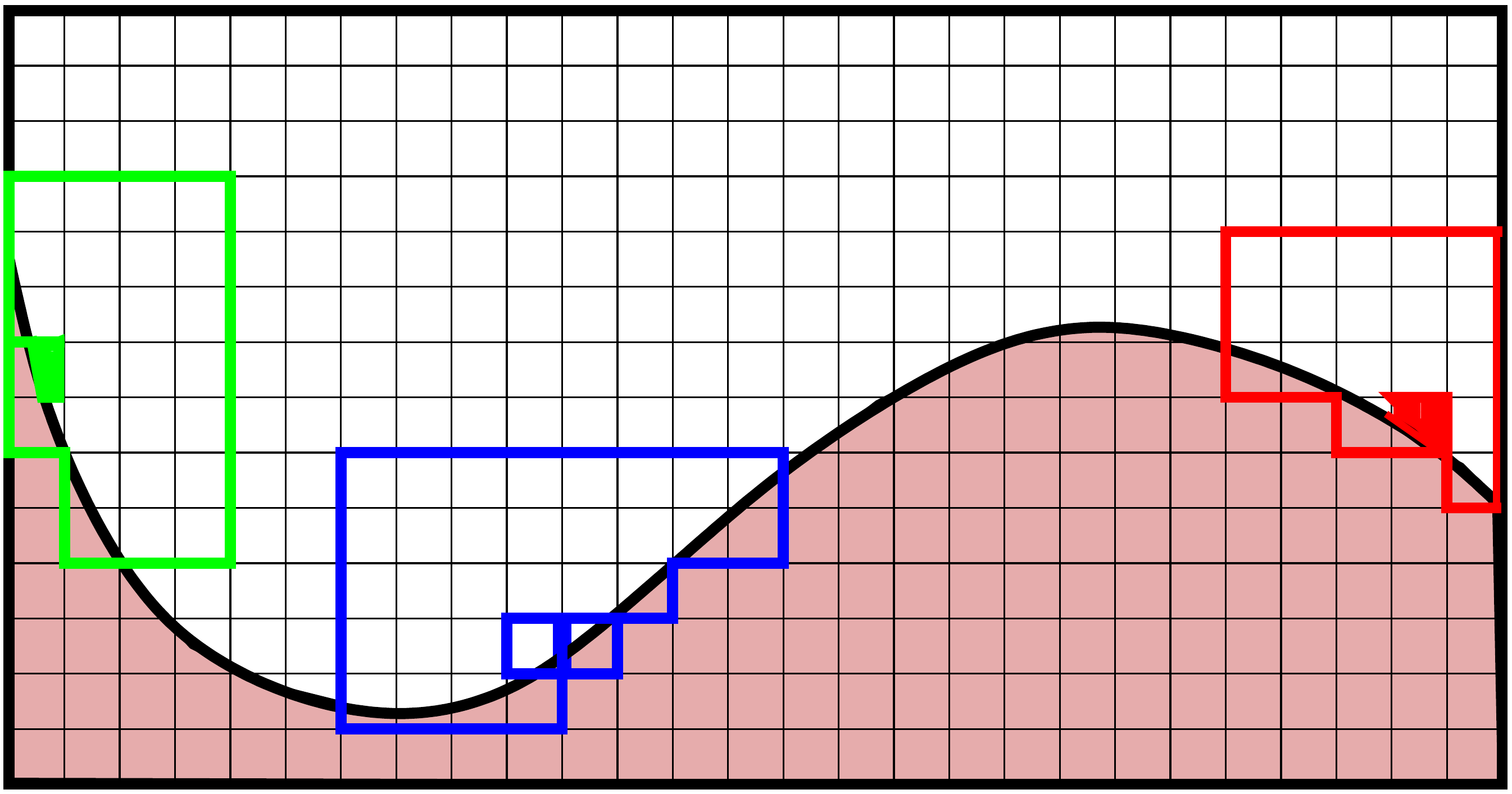}
\caption{Neighbors of a grid-aligned face (in blue) and two embedded boundary faces (whose corresponding cut cells are highlighted in red and green).  In this diagram, the path radius $R_n$ is 3.}
\end{subfigure}
\hspace{0.03\textwidth}
\begin{subfigure}{0.5\textwidth}
\includegraphics[width=\textwidth]{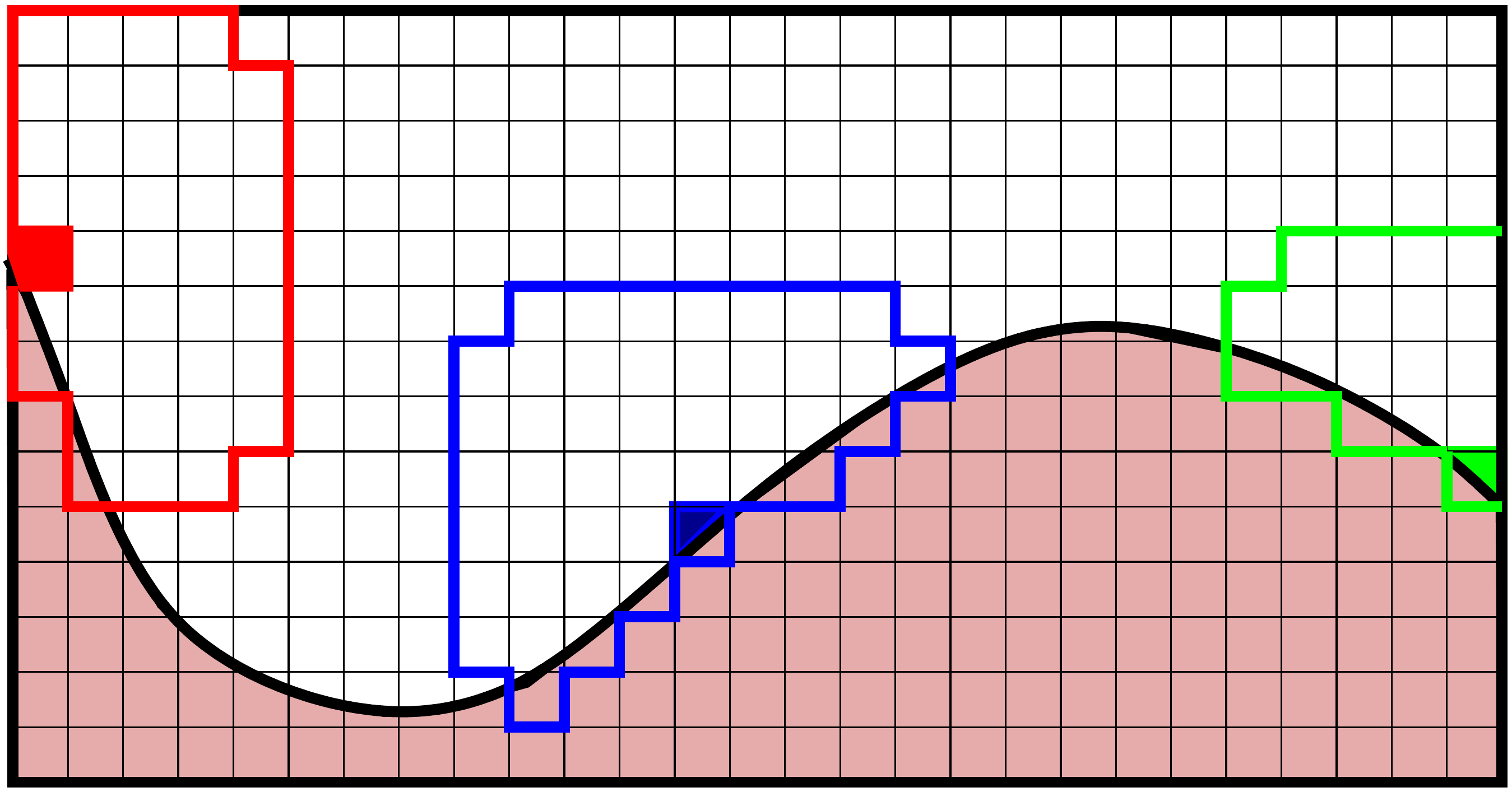}
\caption{Volumes involved in the Laplacian stencils for three volumes (highlighted in blue, green, and red).}
\end{subfigure}
\caption{Neighbors used to construct stencils (flux stencils for faces and Laplacian stencils for volumes).}
\label{fig:stencil ranges}
\end{figure}

We take the neighbors of a volume $V$ to be the volumes in the physical domain that are $R_n$ cells away from $V$.  We call $R_n$ the path radius.  If a neighboring volume touchs a boundary face, then that boundary face is also regarded as a neighbor of $V$.  The neighbors of a face $f$ are the neighbors of its adjacent volumes.  See Figure \ref{fig:stencil ranges}.

For simplicity, our method uses one path radius for all volumes in the domain.  To pin down the interpolant for a boundary face, we need a large path radius.  In our results, we set $R_n = 2$ for the second-order method, and $R_n=3$ for the fourth-order method.  

As explained at the end of Section \ref{sec: polynomial interpolant}, we want a weighting that emphasizes nearest neighbors, and gives little weight to far away neighbors.  Our results were constructed with the weighting matrix $W$, with entries
\begin{equation}
\nonumber
W_{ii} =  
  \begin{cases} 
    1 & \frac{|\xb_i - \xz|}{h}  < \half \\ 
    (2 \frac{|\xb_i - \xz|}{h})^{-5} & \frac{|\xb_i - \xz|}{h} \ge \half 
  \end{cases} 
.
\end{equation}
If the $i$th row of $A$ (equation \refEq{def of A}) corresponds to a volume, then $\xb_i$ is the cell center.  If it corresponds to a grid-aligned boundary face, then $\xb_i$ is the face center.  Finally, if it corresponds to an embedded boundary face, then $\xb_i$ is the cell center of the cut cell.

\section{Results}
\labelSec{results}
%{\tt tests.tex} \\

In this section, we apply the method to solve Poisson's equation on several
different geometries, and demonstrate that the method achieves second and fourth
order accuracy and produces a stable Laplacian operator.  

\subsection{Approximate Moments}
The results presented in this paper were all obtained using $O(h^{6+D-1})$ accurate face moments and $O(h^{6+D})$ accurate volume moments.  (Note that for $D=2$, the grid-aligned face moments are one-dimensional and are calculated exactly.)

\subsection{Solver}
We use the algebraic multigrid (AMG) method in the PETSc solver
framework for our linear equation solves \cite{petsc-web-page,
petsc-user-ref, petsc-efficient}.  PETSc provides interfaces to
several third party AMG solver and has a build in AMG solver, GAMG.
We use a GMRES solver preconditioned with GAMG, which implements a
smoothed aggregation AMG method \cite{Adams-03a}.  To compute
eigenvalues we use SLEPc, an extension of PETSc
\cite{Hernandez:2003:SSL,slepc-users-manual}.  The eigenvalues are
found using Davidson methods \cite{Romero:2014:PID}.

\subsection{Computing Convergence Rates}
To compute convergence rates, we prescribe an exact solution and compare the computed quantities to the exact quantities.  Convergence rates are calculated for the face fluxes, the operator truncation error, and the solution error.  For a face $f$, the error in the face flux is
	\begin{equation}
	\nonumber
		e^{\text{flux}}_f = \frac{1}{A_f}s_f^T\Phi - \frac{1}{A_f}\int_f \nabla \phi \cdot n d\area ,
	\end{equation} 
where $s_f$ is the flux stencil (equation \refEq{stencil equation solution}) and $\Phi$ is defined in equation \refEq{def of phi vec}.  Here, we have written $s_f$ instead of $s$ to highlight that $s$ depends on $f$.  The flux error is computed for all of the faces, including boundary faces.  The (kappa-weighted) truncation error for a volume $\vol_v$ is 
	\begin{equation}
	\nonumber
		e^{\text{lapl}}_v = \frac{\kappa_v}{|\vol_v|}\sum_{f(v)}s^T_{f(v)}\Phi - \frac{\kappa}{|\vol_v|}\int_{\vol_v} \Delta \phi d\vol .
	\end{equation}
Finally, the solution error for $\vol_v$ is the difference between the computed volume-average of $\phi$ and the exact volume-average.

The $L_\infty$ error is the maximum absolute error over the faces or the cells.  For the flux, the $L_1$ error is defined as
	\begin{equation}
	\nonumber
		e^{\text{flux}}_{L_1} = \frac{1}{\sum_{f}\alpha_f}\sum_{f} |e^{\text{flux}}_f| \alpha_f.
	\end{equation}
The $L_1$ truncation error is
	\begin{equation}
	\nonumber
		e^{\text{lapl}}_{L_1} = \sum_{\vol_v} |e^{\text{lapl}}_v| |\vol_v|.
	\end{equation}
Similarily for the solution error.

In each of 2D tests, the exact solution is
	\begin{equation}
	\nonumber
		\phi(x,y) = \sin(2\pi(x-\sqrt{2}/2))\sin(2\pi(x-\sqrt{3}/2)) .
	\end{equation}
The exact solution for the 3D example is given in Section \ref{sec: 3D example}.

%From interpolation theory, the flux will be $O(h^Q)$ locally for a Qth order method.  As proven in (reference to Hans' paper here), the truncation error will be $O(h^{Q-1})$ in cells that use boundary information to generate the stencil.  In the other cells, the truncation error will be $O(h^Q)$ because of symmetry.  Thus, the $L_{\infty}$ truncation error is $O(h^{Q-1})$.  
%
%There are $O(N)$ points near the boundary and $O(N^2)$ that are sufficiently far from the boundary.  So, in the $L_1$ norm, we have
%	\begin{equation}
%		\sum |e_{ij}|\Delta V \approx N^2(h^Q)(h^2) + N(h^{Q-1})(h^2) = h^Q .
%	\end{equation}
%
%The lower order truncation error in cells near the boundary does not reduce the local order of accuracy in the solution.  Using potential theory arguments (see Hans' paper), it can be shown that the error is $O(h^{Q+1})$ near the boundary when Dirichlet boundary conditions are used and $O(h^Q)$ when Neumann boundary conditions.  In the interior, the solution error is $O(h^Q)$, so the overall order of accuracy is $O(h^Q)$ in all norms. 

\subsection{Circle}
Our first test geometry is the circle 
	\begin{equation}
	\nonumber
		(x-0.5)^2 + (y-0.5)^2 = 0.25^2.
	\end{equation}
We apply the method to the region outside of this circle.    

\begin{figure}
\begin{subfigure}{0.5\textwidth}
\includegraphics[width = \textwidth]{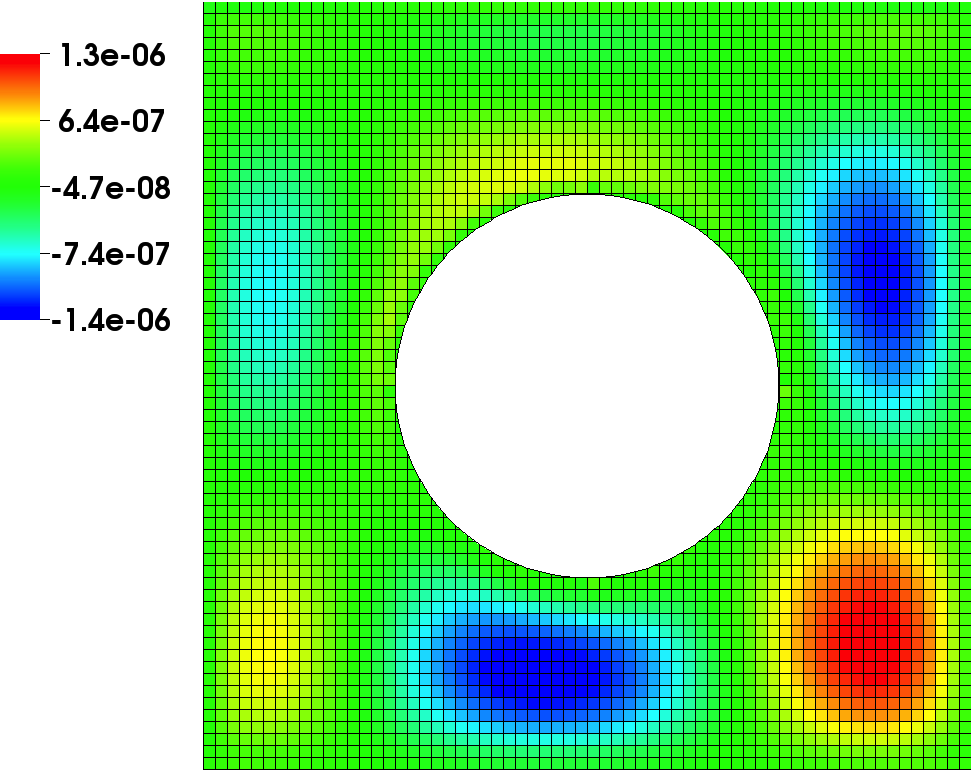}
\caption{Dirichlet boundary conditions on the circle.}
\end{subfigure}
\hspace{0.03\textwidth}
\begin{subfigure}{0.5\textwidth}
\includegraphics[width = \textwidth]{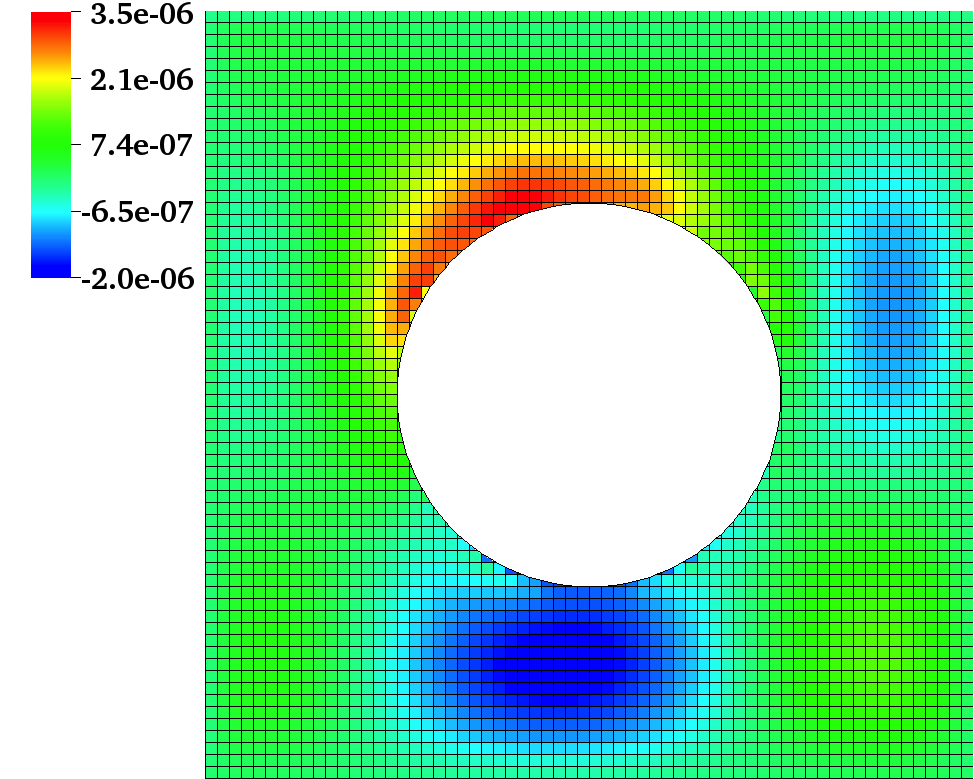}
\caption{Neumann boundary conditions on the circle.}
\end{subfigure}
\caption{Solution error for fourth order method applied to area outside of circle with center $(0.5, 0.5)$ and radius $0.25$.  Dirichlet boundary conditions are prescribed on the square.  The meshwidth of the Cartesian grid is $h = 1/64$.}
\label{fig: circle solution error}
\end{figure}

For this geometry, we simulate three cases:
	\begin{enumerate}
		\item 2nd order method with Dirichlet boundary conditions on the boundaries,
		\item 4th order method with Dirichlet boundary conditions on the boundaries,
		\item 4th order method with Dirichlet boundary conditions on the square and Neumann boundary conditions on the circle.
	\end{enumerate}
Figure \ref{fig: circle solution error} shows the solution error for test cases 2 and 3.  Convergence rates are given in Table \ref{table: circle convergence rates} for the three cases.  

%As expected, the flux and the solution error converge at second or fourth order.  In the $L_1$ norm, the truncation error is second or fourth order as shown in the previous analysis.  
%However, the truncation error in the $L_\infty$ norm does not converge at the expected orders of one or three.  This is because the the truncation error contains a factor of $1/\kappa$, and cut cell volume fractions change rapidly from grid to grid.  Moreover, this example contains small cells: the minimum volume fraction is 1e-2 on the grid with meshwidth $h = 1/64$ and 2.5e-3 on the grid with $h=1/128$.  Thus, the factor of $1/\kappa$ is large, and significantly affects the convergence rate.  If we look at the $\kappa$-weighted truncation error (that is, the truncation error times $\kappa$), we see the expected convergence rates in $L_\infty$ (see subtable \ref{table: circle Linf kappa trunc}).  Because of this dependence on $\kappa$, we only show the $\kappa$-weighted truncation error for the other geometries.

\begin{figure}
\begin{subfigure}{\textwidth}
\includegraphics[width = \textwidth]{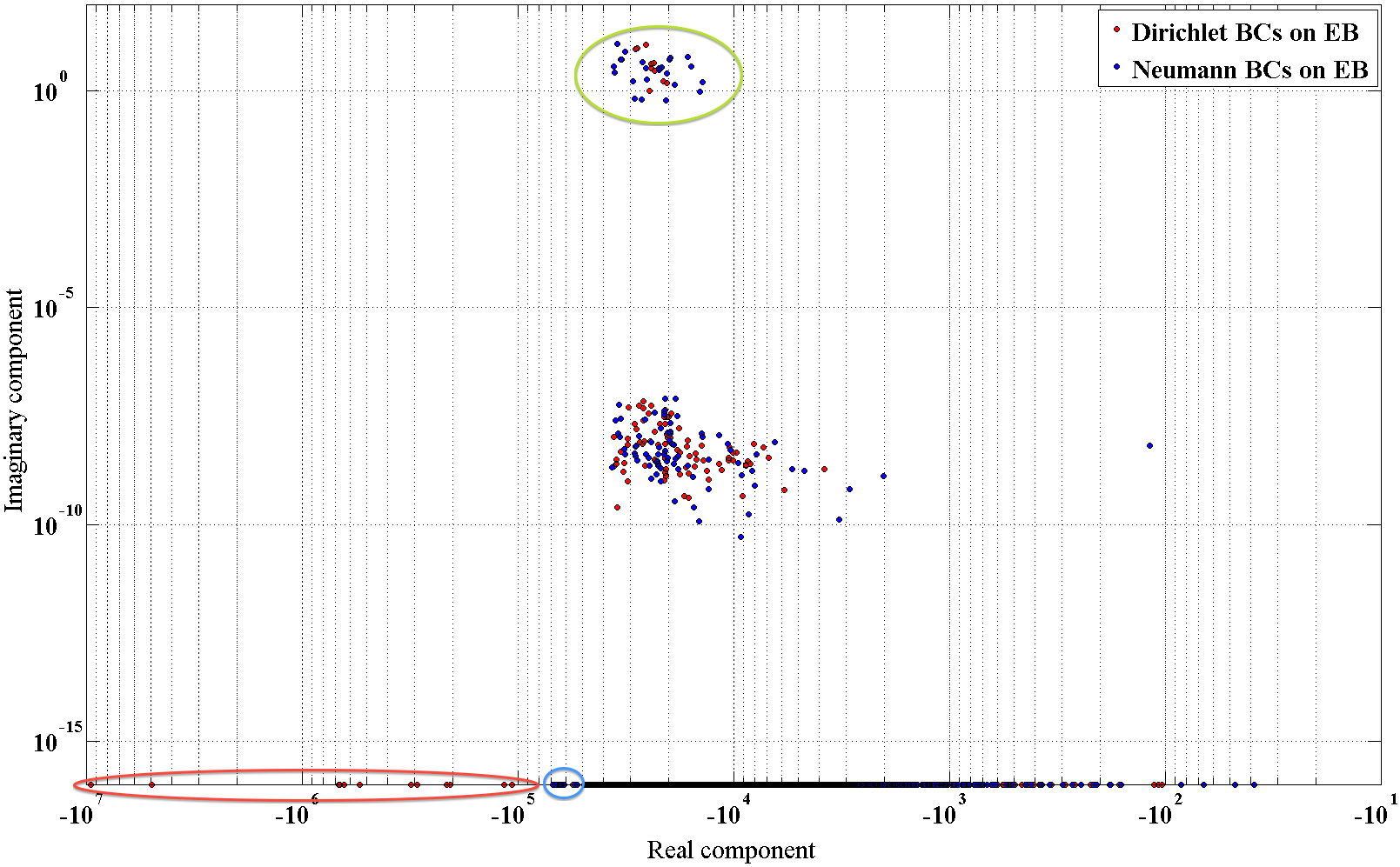}
\caption{Spectrum for fourth order method applied to domain outside of circle wih center $(0.5, 0.5)$ and radius $0.25$.  Dirichlet boundary conditions are prescribed on the square.  The operator is not weighted by volume fraction.  Only the eigenvalues with positive imaginary components are plotted because eigenvalues come in complex conjugate pairs.  Independent of the circle boundary condition type, the operator has a few eigenvalues with large imaginary component (indicated by green circle).  It also has large negative eigenvalues (indicated by the red circle for Dirchlet boundary conditions and the blue circle for Neumann boundary conditions) because small cells are present.  The Cartesian grid is $64\times 64$.}
\label{fig: circle spectrum}
\end{subfigure}
\hspace{0.01\textwidth}
\begin{subfigure}{\textwidth}
\begin{center}
\begin{tabular}{| c | c | c | c |}
\hline
Test Case & $\lambda^{*}$ & $\lambda_{\text{max}}$ & $\lambda_{\text{min}}$ \\ \hline
$x_0=(0.5,0.5)$, $r = 0.255$, Dirichlet BCs on EB &-3.6e4 + 28$i$ & -106 & -1.0e6\\ \hline
$x_0=(0.501,0.501)$, $r = 0.25$, Dirichlet BCs on EB & -3.1e4 + 5.0$i$ & -103 & -2.5e6 \\ \hline
%$x_0=(0.51,0.5)$, $r = 0.25$, Dirichlet BCs on EB & -3.8e6 + 2.5e6$i$ & -101 & -1.7e7\\ \hline
$x_0=(0.5,0.5)$, $r = 0.255$, Neumann BCs on EB & -3.0e4 + 24$i$ &-39 &-9.1e4  \\ \hline
$x_0=(0.501,0.501)$, $r = 0.25$, Neumann BCs on EB & -3.5e4 + 17$i$& -38 &-1.3e5 \\ \hline
$x_0=(0.51,0.5)$, $r = 0.25$, Neumann BCs on EB & -2.0e4 + 10$i$ & -38 & -4.0e5\\
\hline
\end{tabular}
\end{center}
\caption{Summary of spectra information for the fourth order method
  applied to the circle geometries.  $x_0$ is the center of the
  circle, $r$ is the radius.  The three columns are: the eigenvalue
  with the largest imaginary component ($\lambda^{*}$), the eigenvalue
  with the largest real component ($\lambda_{\text{max}}$), and the
  eigenvalue with the smallest real component ($\lambda_{\text{min}}$).  Note that the eigenvalues come in complex conjugate pairs, and only the eigenvalue with a positive imaginary component is given in the first column.}
\end{subfigure}
\caption{}
\label{fig: geometric perturbation spectra summary}
\end{figure}

Figure \ref{fig: circle spectrum} shows the spectra for the Laplacian operator in test cases 2 and 3.  We have not weighted the operator by the volume fraction.  In both cases the eigenvalues of the operator lie in the left half plane.  There are a few eigenvalues with large negative real part.  The eigenvectors corresponding to these eigenvalues are concentrated in the small cut cells.  There is also a cluster of eigenvalues whose imaginary parts are relatively large.  If this operator is combined with a temporal discretization like the trapezoidal rule to solve the heat equation, these eigenvalues would introduce oscillations in the solution; however, these oscillations are quickly damped out because the eigenvalues also have large negative real parts.

Without weights, the operator has large positive eigenvalues and is
unstable (data not shown).  The eigenvectors for the positive
eigenvalues are supported on the small cells.  The operator also has a
pair of eigenvalues with significantly large imaginary parts, with eigenvectors that
couple the smallest cells.

\begin{figure}
\includegraphics[width = \textwidth]{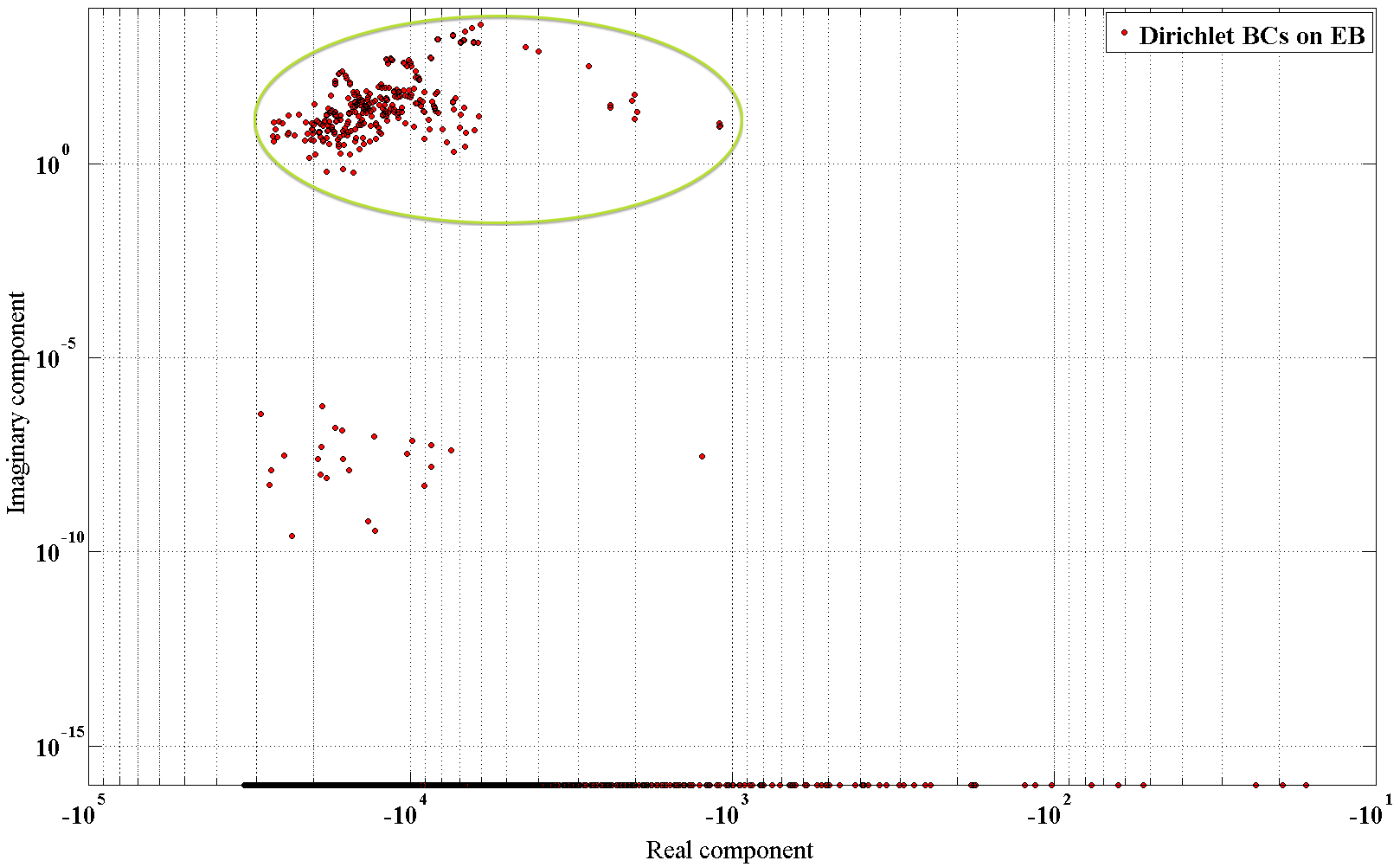}
\caption{Spectrum for operator generated using the method in
  \cite{schwartzETAL:2006}.  The operator has been weighted by volume
  fraction.  The domain is the region outside of the circle wih center $(0.5, 0.5)$ and radius $0.25$.  Dirichlet boundary conditions are prescribed on the boundaries.  Like the operator for the weighted least squares method, this operator has a cluster of eigenvalues with non-negligible imaginary component (indicated by the green circle).  The Cartesian grid is $64\times 64$.}
\label{fig: ebamrpoisson spectrum}
\end{figure}

As a point of reference, we have plotted the spectrum for the second-order operator generated using the method in \cite{schwartzETAL:2006}, (see Figure \ref{fig: ebamrpoisson spectrum}).  This operator is weighted by the volume fraction.  Like our fourth-order operator, this operator also has eigenvalues with non-negligible imaginary components. 

%Finally, Figure \ref{fig: unweighted least squares spectrum} shows the spectrum for the operator generated with standard least squares.  The operator has not been multiplied by the matrix of volume fractions.  Without weights, the operator has large positive eigenvalues and is unstable.  

% begin table
\begin{table}
\begin{subtable}{\textwidth}
\begin{center}
  \begin{tabular}{| c | c | c | c | c | c |}
    \hline
    Test & $e_{32}$ & Order & $e_{64}$ & Order & $e_{128}$\\ \hline
    2nd order, Dirichlet EB &  5.01 & 1.03  &   2.44 & 1.72 & 7.43e-01  \\ \hline
    4th order, Dirichlet EB &   9.36e-02 & 3.96  &   6.00e-03 & 2.62 & 9.72e-04   \\ \hline
    4th order, Neumann EB &  5.30e-02 & 2.00   &  1.33e-02 & 2.83 & 1.87e-03  \\
    \hline
  \end{tabular}
\end{center}
\caption{$L_\infty$ convergence rates for $\kappa$-weighted truncation error.}
\label{table: circle Linf kappa trunc}
\end{subtable}
\\
\begin{subtable}{\textwidth}
\begin{center}
  \begin{tabular}{| c | c | c | c | c | c |}
    \hline
    Test & $e_{32}$ & Order & $e_{64}$ & Order & $e_{128}$\\ \hline
    2nd order, Dirichlet EB & 2.39e-01 & 1.91  &   6.39e-02 & 1.96 & 1.64e-02  \\ \hline
    4th order, Dirichlet EB &   4.89e-03 & 3.97  &   3.11e-04 & 3.89 & 2.10e-05  \\ \hline
    4th order, Neumann EB &  5.37e-03 & 3.87 &  3.68e-04 & 3.80 & 2.64e-05  \\
    \hline
  \end{tabular}
\end{center}
\caption{$L_1$ convergence rates for $\kappa$-weighted truncation error.}
\label{table: circle L1 kappa trunc}
\end{subtable}
\\
\begin{subtable}{\textwidth}
\begin{center}
  \begin{tabular}{| c | c | c | c | c | c |}
    \hline
    Test & $e_{32}$ & Order & $e_{64}$ & Order & $e_{128}$\\ \hline
    2nd order, Dirichlet EB &   1.36e-03 & 1.90  &   3.65e-04 & 1.92 & 9.64e-05  \\ \hline
    4th order, Dirichlet EB &   2.56e-05 & 4.17 & 1.43e-06 & 3.86 & 9.80e-08   \\ \hline
    4th order, Neumann EB &  5.24e-05 & 3.89   &  3.52e-06 & 3.94 & 2.29e-07  \\
    \hline
  \end{tabular}
\end{center}
\caption{$L_\infty$ convergence rates for solution error.}
\label{table: circle Linf soln error}
\end{subtable}
\\
\begin{subtable}{\textwidth}
\begin{center}
  \begin{tabular}{| c | c | c | c | c | c |}
    \hline
    Test & $e_{32}$ & Order & $e_{64}$ & Order & $e_{128}$\\ \hline
    2nd order, Dirichlet EB &  4.06e-04 & 2.04  &  9.90e-05 & 1.97 & 2.52e-05 \\ \hline
    4th order, Dirichlet EB & 5.70e-06 & 3.91  & 3.78e-07 & 3.91 & 2.52e-08 \\ \hline
    4th order, Neumann EB &  1.07e-05 & 4.05  & 6.48e-07 & 3.96 & 4.16e-08\\
    \hline
  \end{tabular}
\end{center}
\caption{$L_1$ convergence rates for solution error.}
\label{table: circle L1 soln error}
\end{subtable}
\caption{Convergence rates for method applied to area outside of circle with center $(0.5, 0.5)$ and radius $0.25$.  $e_N$ is the error for a grid with meshwidth $h=1/N$.  Dirichlet boundary conditions are prescribed on the square.}
\label{table: circle convergence rates}
\end{table}
% end table

%\begin{figure}
%\includegraphics[width = \textwidth]{testFigs/unweightedLS.pdf}
%\caption{Spectrum for fourth order \emph{unweighted} least squares method applied to domain outside of circle wih center $(0.5, 0.5)$ and radius $0.25$.  Dirichlet boundary conditions are prescribed on the domain boundary.  The operator is not weighted by volume fraction.  The Cartesian grid is $64\times 64$.}
%\label{fig: unweighted least squares spectrum}
%\end{figure}

%%%%%%%%%%%%%%%%%%%%%%%%%
\subsubsection{Geometric Perturbations}
Perturbations in the geometry can produce small cells.  Our next set of tests demonstrates that the method is robust under small changes in the geometry.  In particular, we perturb the radius and center of the circle in the last example, and apply the fourth order method to solve Poisson's equation on the perturbed domains.  Table \ref{table: smallest volume fractions for geometric perturbations} lists the smallest volume fractions for the original circle example and three perturbations of this circle on the three grids used in the study.  Note that perturbing the circle's center from $x_0 = (0.5, 0.5)$ to $x_0 = (0.51, 0.5)$ changes the smallest volume fraction by a magnitude!

Despite these differences in the geometry, the $\kappa$-weighted
truncation error and solution error hardly change.  Tables \ref{table:
  geometric perturbation dirichlet convergence rates} and \ref{table:
  geometric perturbation neumann convergence rates} display the
convergence rates in the case of Dirichlet boundary conditions and
Neumann boundary conditions on the embedded boundary, respectively.

The spectrum for each perturbation is similar to the spectrum in the
previous example.  The eigenvalues lie in the left half-plane for each
of the perturbations.  Also, there are a few eigenvalues with large
negative real part corresponding to the small cells and a cluster of
eigenvalues with non-negligible imaginary part.  We summarize the main features of each spectrum in Figure \ref{fig: geometric perturbation spectra summary}.    

%There is one significant difference between the spectra: the circle with $x_0 = (0.51, 0.5)$ and $r = 0.25$ has two (complex-conjugate) eigenvalues with large imaginary components.  For this particular geometry, the top of the circle slices the grid and generates two neighboring small thin cells.  These cells are connected by a face with area fraction $1/300$.  The eigenvectors corresponding to those eigenvalues with large imaginary components are concentrated in the two small cells.  Such a situation with two thin cells does not occur in the other three cases.

\begin{table}
\begin{center}
\begin{tabular}{|c | c | c | c |}
\hline
Geometry & N=32 & N=64 & N=128 \\ \hline
$x_0 = (0.5, 0.5)$, $r= 0.25$ &  4.5e-3 &  1.0e-2 & 2.5e-4\\ \hline
$x_0 = (0.5, 0.5)$, $r= 0.255$ & 1.7e-2 & 6.8e-3 & 9.5e-5 \\ \hline
$x_0 = (0.501, 0.501)$, $r= 0.25$ & 4.0e-4 & 1.6e-3 & 7.2e-6 \\ \hline
$x_0 = (0.51, 0.5)$, $r= 0.25$ & 3.6e-4 & 2.3e-4 & 4.5e-5 \\
\hline
\end{tabular}
\end{center}
\caption{Smallest volume fraction on a $N$x$N$ Cartesian grid for the domain outside of circle with center $x_0$ and radius $r$. }
\label{table: smallest volume fractions for geometric perturbations}
\end{table}

\begin{table}
\begin{subtable}{\textwidth}
\begin{center}
  \begin{tabular}{| c | c | c | c | c | c |}
    \hline
    Test & $e_{32}$ & Order & $e_{64}$ & Order & $e_{128}$\\ \hline
    $x_0 = (0.5, 0.5)$, $r= 0.255$ &    7.24e-02 & 3.59  &   6.00e-03 & 3.10 & 7.01e-04  \\ \hline
    $x_0 = (0.501, 0.501)$, $r=0.25$ &   5.20e-02 & 3.12 & 6.00e-03   &  3.10 & 7.01e-04 \\ \hline
    $x_0 = (0.51, 0.5) $, $r= 0.25$ &   7.41e-02 & 3.46 & 6.71e-03  &  2.39 & 1.28e-03 \\
    \hline
  \end{tabular}
\end{center}
\caption{$L_\infty$ convergence rates for $\kappa$-weighted truncation error}
\end{subtable}
\\

\begin{subtable}{\textwidth}
\begin{center}
  \begin{tabular}{| c | c | c | c | c | c |}
    \hline
    Test & $e_{32}$ & Order & $e_{64}$ & Order & $e_{128}$\\ \hline
    $x_0 = (0.5, 0.5)$, $r= 0.255$ &   4.93e-03 & 3.92  &  3.26e-04 & 3.91 & 2.17e-05 \\ \hline
    $x_0 = (0.501, 0.501)$, $r=0.25$ &  4.87e-03 & 3.95 & 3.14e-04  & 3.87 & 2.14e-05\\ \hline
    $x_0 = (0.51, 0.5) $, $r= 0.25$ &  5.06e-03 & 3.98 & 3.21e-04 & 3.91 & 2.14e-05 \\
    \hline
  \end{tabular}
\end{center}
\caption{$L_1$ convergence rates for $\kappa$-weighted truncation error}
\end{subtable}
\\
\begin{subtable}{\textwidth}
\begin{center}
  \begin{tabular}{| c | c | c | c | c | c |}
    \hline
    Test & $e_{32}$ & Order & $e_{64}$ & Order & $e_{128}$\\ \hline
    $x_0 = (0.5, 0.5)$, $r= 0.255$ &   2.01e-05 & 3.85 & 1.40e-06   &  3.88 & 9.48e-08  \\ \hline
    $x_0 = (0.501, 0.501)$, $r=0.25$ &  3.38e-05 & 4.55 & 1.45e-06   &  3.88 & 9.84e-08 \\ \hline
    $x_0 = (0.51, 0.5) $, $r= 0.25$ &  3.02e-05 & 3.34 & 2.97e-06  &  4.91 & 9.88e-08  \\
    \hline
  \end{tabular}
\end{center}
\caption{$L_\infty$ convergence rates for solution error}
\end{subtable}
\\
\begin{subtable}{\textwidth}
\begin{center}
  \begin{tabular}{| c | c | c | c | c | c |}
    \hline
    Test & $e_{32}$ & Order & $e_{64}$ & Order & $e_{128}$\\ \hline
    $(x_0, y_0) = (0.5, 0.5) $, $r=0.255$ &  5.56e-06 & 3.90  &  3.73e-07 & 3.92 & 2.47e-08 \\ \hline
    $(x_0, y_0) = (0.501, 0.501)$, $r=0.25$ &   6.06e-06 & 4.00   &  3.79e-07 & 3.91 & 2.51e-08\\ \hline
    $(x_0, y_0) = (0.51, 0.5)$, $r=0.25$ &  5.94e-06 & 3.93 & 3.90e-07 & 3.96 & 2.50e-08\\
    \hline
  \end{tabular}
\end{center}
\caption{$L_1$ convergence rates for solution error}
\end{subtable}
\caption{Convergence rates for 4th order method applied to domain outside of the circle $(x-x_0)^2 + (y-y_0)^2 = r^2$.  Dirichlet boundary conditions are prescribed on the boundaries.}
\label{table: geometric perturbation dirichlet convergence rates}
\end{table}

\begin{table}
\begin{subtable}{\textwidth}
\begin{center}
  \begin{tabular}{| c | c | c | c | c | c |}
    \hline
    Test & $e_{32}$ & Order & $e_{64}$ & Order & $e_{128}$\\ \hline
    $x_0 = (0.5, 0.5)$, $r= 0.255$ &    5.19e-02 & 2.29 & 1.06e-02   &  2.48 & 1.90e-03 \\ \hline
    $x_0 = (0.501, 0.501)$, $r=0.25$ &  5.17e-02 & 2.03 & 1.27e-02  & 2.69 & 1.96e-03 \\ \hline
    $x_0 = (0.51, 0.5) $, $r= 0.25$ &    5.53e-02 & 2.09 & 1.30e-02  &  2.63 & 2.10e-03  \\
    \hline
  \end{tabular}
\end{center}
\caption{$L_\infty$ convergence rates for $\kappa$-weighted truncation error}
\end{subtable}
\\

\begin{subtable}{\textwidth}
\begin{center}
  \begin{tabular}{| c | c | c | c | c | c |}
    \hline
    Test & $e_{32}$ & Order & $e_{64}$ & Order & $e_{128}$\\ \hline
    $x_0 = (0.5, 0.5)$, $r= 0.255$ &    5.41e-03 & 3.77 & 3.95e-04 &  3.90 & 2.64e-05 \\ \hline
    $x_0 = (0.501, 0.501)$, $r=0.25$ &   5.42e-03 & 3.86 & 3.72e-04 &  3.85 & 2.581e-05\\ \hline
    $x_0 = (0.51, 0.5) $, $r= 0.25$ &   5.40e-03 & 3.82 & 3.83e-04 & 3.91 & 2.542e-05 \\
    \hline
  \end{tabular}
\end{center}
\caption{$L_1$ convergence rates for $\kappa$-weighted truncation error}
\end{subtable}
\\
\begin{subtable}{\textwidth}
\begin{center}
  \begin{tabular}{| c | c | c | c | c | c |}
    \hline
    Test & $e_{32}$ & Order & $e_{64}$ & Order & $e_{128}$\\ \hline
    $x_0 = (0.5, 0.5)$, $r= 0.255$ &   5.41e-05 & 3.95 & 3.49e-06   &  3.94 & 2.27e-07  \\ \hline
    $x_0 = (0.501, 0.501)$, $r=0.25$ &  5.32e-05 & 3.91 & 3.53e-06   &  3.95 & 2.29e-07 \\ \hline
    $x_0 = (0.51, 0.5) $, $r= 0.25$ &  5.52e-05 & 3.95 & 3.58e-06  &  3.94 & 2.33e-07  \\
    \hline
  \end{tabular}
\end{center}
\caption{$L_\infty$ convergence rates for solution error}
\end{subtable}
\\
\begin{subtable}{\textwidth}
\begin{center}
  \begin{tabular}{| c | c | c | c | c | c |}
    \hline
    Test & $e_{32}$ & Order & $e_{64}$ & Order & $e_{128}$\\ \hline
    $x_0 = (0.5, 0.5) $, $r=0.255$ &  1.09e-05 & 4.08  & 6.45e-07 & 3.97 & 4.12e-08 \\ \hline
    $x_0 = (0.501, 0.501)$, $r=0.25$ &  1.08e-05 & 4.06    &   6.48e-07 & 3.96 & 4.17e-08 \\ \hline
    $x_0 = (0.51, 0.5)$, $r=0.25$ &   1.11e-05 & 4.08  &  6.55e-07 & 3.96 & 4.21e-08 \\
    \hline
  \end{tabular}
\end{center}
\caption{$L_1$ convergence rates for solution error}
\end{subtable}
\caption{Convergence rates for 4th order method applied to domain outside of the circle $(x-x_0)^2 + (y-y_0)^2 = r^2$.  Dirichlet boundary conditions are prescribed on the square, and Neumann boundary conditions are prescribed on the circle.}
\label{table: geometric perturbation neumann convergence rates}
\end{table}

%%%%%%%%%%%%%%%%%%%%%%%%%
\subsection{Other Geometries}

In the next few sections, we apply the fourth order method on
different geometries.  In each case, the operator spectrum looks
similar to the circle example: there are a few eigenvalues with large
negative part that correspond to small cells and a cluster of
eigenvalues with non-negligible imaginary part.  As a result, we
summarize the main features of the spectra in Table \ref{table:
  spectra summary}.
% need this line because LaTex can only support a small number of floats without text
\clearpage

%%%%%% Sine curve %%%%%%%%%

\subsubsection{Trignometric Curve}
\begin{figure}
\includegraphics[width = \textwidth]{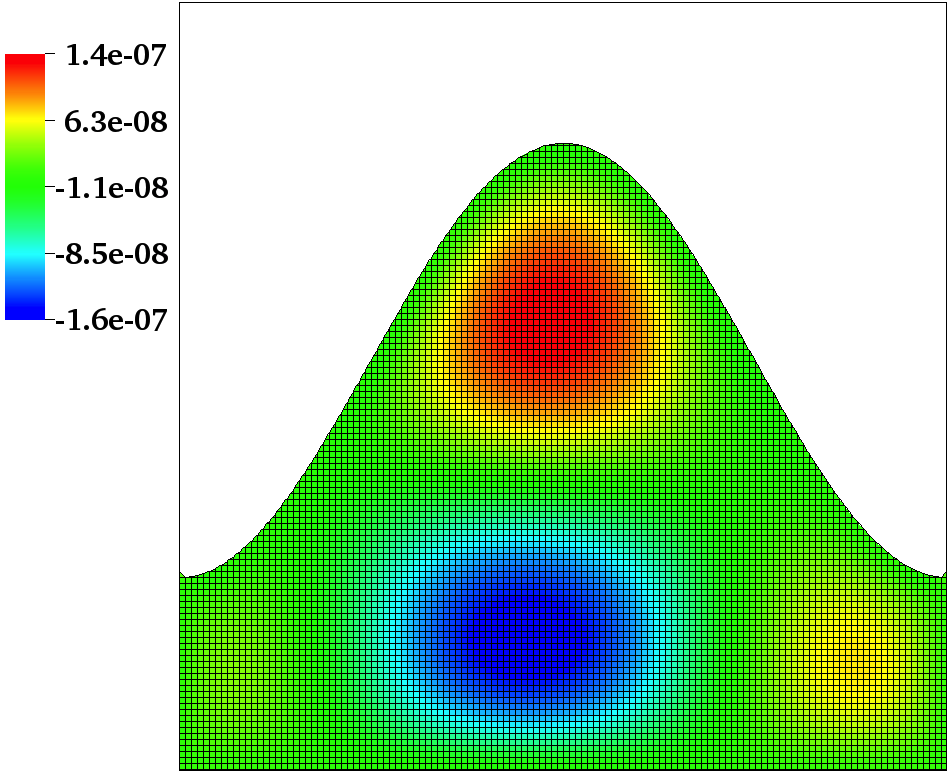}
\caption{Solution error from fourth order method applied to area underneath the curve  $y = 0.25 + \frac{\sqrt{2}}{2}(1-\cos(2\pi x)) $.  Dirichlet boundary conditions are prescribed on the boundaries.  The Cartesian grid is $128\times128$.}
\label{fig: sine curve solution error}
\end{figure}

Our next domain is the region underneath the curve
	\begin{equation}
		y = 0.25 + \frac{\sqrt{2}}{2}(1-\cos(2\pi x)) .
	\end{equation}
We apply the fourth order method to solve Poisson's equation on this geometry.  Dirichlet boundary conditions are prescribed on the boundaries.  Figure \ref{fig: sine curve solution error} shows the solution error on the grid with meshwidth $h = 1/128$, and Table \ref{table: sine curve convergence rates} lists the convergence rates.  It is clear that the method converges at fourth order for this example.

\begin{table}
\begin{center}
  \begin{tabular}{| c | c | c | c | c | c |}
    \hline
    Test & $e_{32}$ & Order & $e_{64}$ & Order & $e_{128}$\\ \hline
    $\kappa$-weighted truncation error, $L_{\infty} $ &  7.82e-02 & 3.38 & 7.53e-03 & 3.06 & 9.03e-04   \\ \hline
    $\kappa$-weighted truncation error, $L_1$ & 3.99e-03 & 3.75 & 2.97e-04 &  3.86 & 2.05e-05  \\ \hline
    Solution error, $L_{\infty} $ &   4.03e-05 & 3.99 & 2.53e-06 & 4.00 & 1.58e-07  \\ \hline
    Solution error, $L_1$ &    1.10e-05 & 3.91 & 7.32e-0 & 3.93 & 4.82e-08   \\ 
    \hline
  \end{tabular}
\end{center}
\caption{Convergence rates for fourth order method applied to area underneath the curve $y = 0.25 + \frac{\sqrt{2}}{2}(1-\cos(2\pi x)) $.  Dirichlet boundary conditions are prescribed on the boundaries.}
\label{table: sine curve convergence rates}
\end{table}

%%%%% Four circles example %%%%%%%

\subsubsection{Four circles}

\begin{figure}
\includegraphics[width = \textwidth]{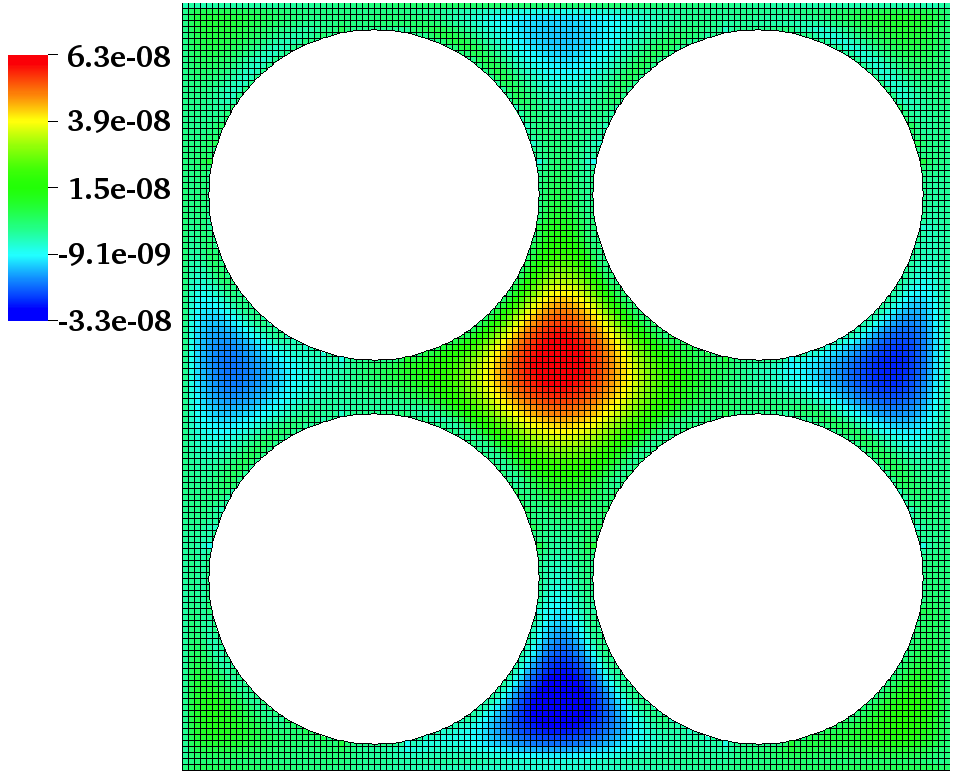}
\caption{Solution error from 4th order applied to the area outside of four circles that are close to the square boundary.  Dirichlet boundary conditions are prescribed on the boundaries.  The Cartesian grid is $128\times 128$.}
\label{fig: four circles solution error}
\end{figure}

Next, we consider the domain outside of the four circles with centers $(0.25, 0.25)$, $(0.75, 0.25)$, $(0.25, 0.75)$, $(0.75, 0.75)$ and all with radius $0.215$.  We apply the fourth order method; Dirichlet boundary conditions are prescribed on the boundaries.  Figure \ref{fig: four circles solution error} shows the solution error on the grid with meshwidth $1/128$.

Unlike the previous examples, we simulate this example on grids with $N = $64, 128, and 256 because the $N = 32$ is too coarse to compute a fourth order flux.   On the $64x64$ grid, each of the circles is only $3$ grid cells away from the square boundary.  Thus, a face close to the boundary has the minimum information to construct a fourth order flux: 3 cells, 1 embedded boundary face, and 1 domain boundary face.  Despite this limited information, the $\kappa$-weighted truncation error and solution error converge at fourth order.  Also, the operator is stable (see Table \ref{table: spectra summary}).  

\begin{table}
\begin{center}
  \begin{tabular}{| c | c | c | c | c | c |}
    \hline
    Test & $e_{64}$ & Order & $e_{128}$ & Order & $e_{256}$\\ \hline
    Truncation error with kappa, $L_{\infty}$ &  3.45e-02 & 5.24 & 9.13e-04 & 3.19 & 1.00e-04\\ \hline
    Truncation error with kappa, $L_1$ & 1.12e-03 & 4.61 & 4.62e-05  & 3.90 & 3.10e-06 \\ \hline
    Solution error, $L _{\infty}$ & 1.04e-06 & 4.04 & 6.30e-08  & 3.94 & 4.11e-09 \\ \hline
    Solution error, $L_1$ &  1.73e-07 & 4.12 & 9.98e-09 & 3.88 & 6.77e-10   \\ 
    \hline
  \end{tabular}
\end{center}
\caption{Convergence rates for fourth order method where the domain is the outside of four circles that are close to the square boundary.}
\label{table: four circle convergence rates}
\end{table}

%%%%%%% SPECTRUM TABLE %%%%%%%%%
\begin{table}
\begin{center}
\begin{tabular}{| c | c | c | c |}
\hline
Test Case & $\lambda^{*}$ & $\lambda_{\text{max}}$ &
$\lambda_{\text{min}}$ \\ \hline
$x_0=(0.5,0.5)$, $r = 0.25$, Dirichlet BCs on EB & -2.5e4 + 11$i$ & -104  & -9.5e6\\ \hline
$x_0=(0.5,0.5)$, $r = 0.255$, Dirichlet BCs on EB &-3.6e4 + 28$i$ & -106 & -1.0e6\\ \hline
$x_0=(0.501,0.501)$, $r = 0.25$, Dirichlet BCs on EB & -3.1e4 + 5.0$i$ & -103 & -2.5e6 \\ \hline
%$x_0=(0.51,0.5)$, $r = 0.25$, Dirichlet BCs on EB & -3.8e6 + 2.5e6$i$ & -101 & -1.7e7\\ \hline
$x_0=(0.5,0.5)$, $r = 0.25$, Neumann BCs on EB & -3.5e4 + 12$i$ & -38  &-6.9e4 \\ \hline
$x_0=(0.5,0.5)$, $r = 0.255$, Neumann BCs on EB & -3.0e4 + 24$i$ &-39 &-9.1e4  \\ \hline
$x_0=(0.501,0.501)$, $r = 0.25$, Neumann BCs on EB & -3.5e4 + 17$i$& -38 &-1.3e5 \\ \hline
$x_0=(0.51,0.5)$, $r = 0.25$, Neumann BCs on EB & -2.0e4 + 10$i$ & -38 & -4.0e5\\ \hline
Other geometries & & & \\ \hline
Four circles, Dirichlet BCs on EB & -4.4e4 + 1.6e3$i$ &  -240 & -6.4e5\\ \hline
Sine curve, Dirichlet BCs on EB & -3.3e4 + 600$i$& -51 & -1.8e6 \\ 
\hline
\end{tabular}
\end{center}
\caption{Summary of spectra information for the 4th order method applied to various geometries.  The three columns are: the eigenvalue
  with the largest imaginary component ($\lambda^{*}$), the eigenvalue
  with the largest real component ($\lambda_{\text{max}}$), and the
  eigenvalue with the smallest real component
  ($\lambda_{\text{min}}$).  Note that the eigenvalues come in complex
  conjugate pairs, and only the eigenvalue with a positive imaginary
  component is given in the first column. The meshwidth of the Cartesian grid is $1/64$.  Dirichlet boundary conditions are prescribed on the square boundary.  }
\label{table: spectra summary}
\end{table}

%%%%%%%%%%%%%%%%%%%%%%%%%
\subsection{Sphere} \label{sec: 3D example}
\begin{figure}
\includegraphics[width = 1.25\textwidth]{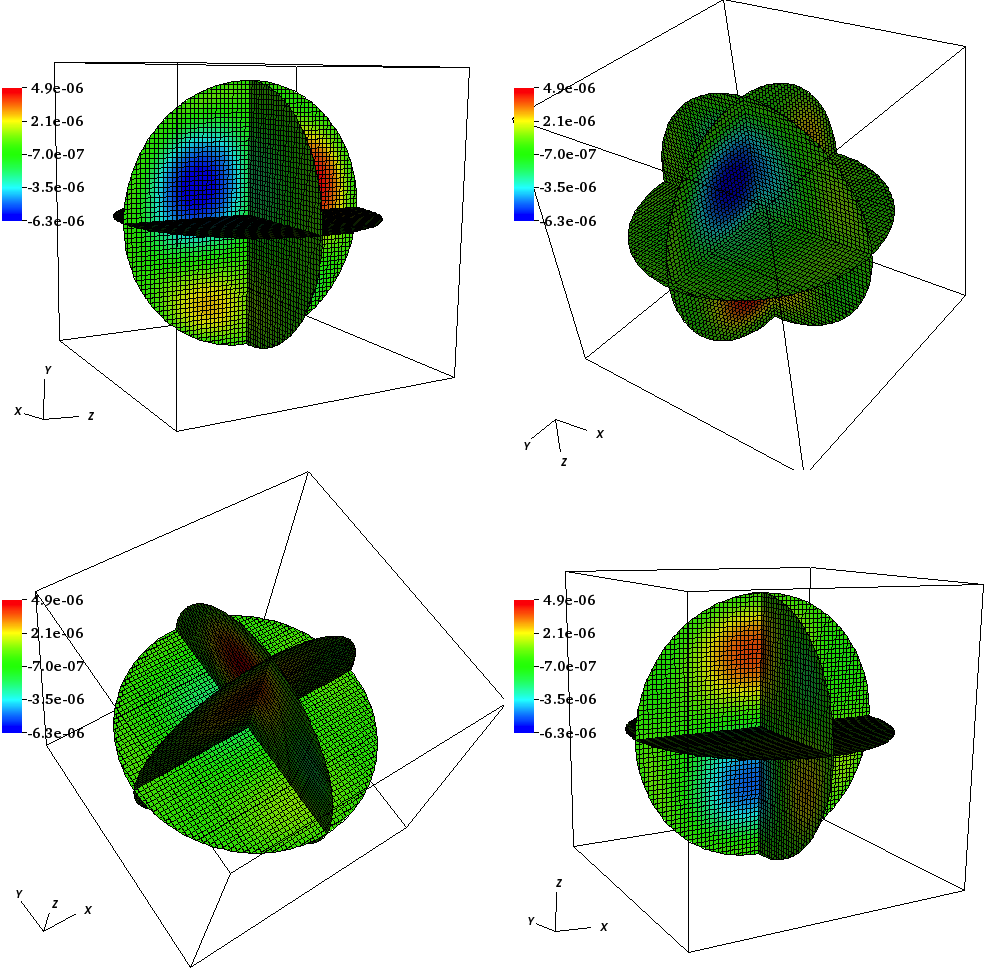}
\caption{Solution error for 4th order method applied to inside of sphere  with center $(0.5, 0.5, 0.5)$ and radius $0.45$.  Each image shows three cross-sections through the center of the sphere.  The triad next to each image shows the orientation of the sphere in the image.}
\label{fig: 3D soln error for sphere}
\end{figure}

In this section, we apply a fourth order method to solve Poisson's equation on the inside of a sphere, with center $(0.5, 0.5, 0.5)$ and radius 0.45.  The exact solution for this test problem is
	\begin{equation}
		\phi(x,y,z) = \sin(2\pi(x-x_0))\sin(2\pi(y-y_0))\sin(2\pi(z-z_0)),
	\end{equation}
where $x_0 = \sqrt{2}/2$, $y_0 = \sqrt{3}/2$, and $z_0 = 0.54321$.  Table \ref{table: 3D soln error for sphere} shows the solution error convergence rates.  It is clear that the method is fourth order.  Figure \ref{fig: 3D soln error for sphere} shows the solution error on a grid with meshwidth $h = 1/128$.  

\begin{table}[p]
\begin{center}
  \begin{tabular}{| c | c | c | c |}
    \hline
    Norm & $e_{32}$ & Order & $e_{64}$ \\ \hline
	$L_\infty$ & 9.16e-05 &    3.87 & 6.27e-06 \\ \hline
	$L_1$ &  1.38e-05 &    3.81 & 9.81e-07\\ \hline
	$L_2$ &  2.17e-05 &    3.81 & 1.54e-06 \\ 
    \hline
  \end{tabular}
\end{center}
\caption{Solution error convergence rates in the $L_\infty$, $L_1$, and $L_2$ norms for 4th order method applied to inside of sphere with center $(0.5, 0.5, 0.5)$ and radius $0.45$.}
\label{table: 3D soln error for sphere}
\end{table}

\section{Conclusions}
\labelSec{conclusions}
%{\tt conclusions.tex} \\

We have presented an algorithm to generate higher-order conservative finite volume discretizations for Poisson's equation on cut cell grids.   The Poisson operator is written in terms of face fluxes, which we approximate using a polynomial interpolant.  The key to the method is to use weighted least squares to generate stable stencils.  In particular, the linear system for the face flux stencil is underdetermined, and we can use weights to pick a stable stencil from the space of solutions.

By applying the method to a variety of geometries, we have demonstrated that the method achieves second and fourth order accuracy.  In each of these examples, we have also shown that the discrete Laplacian operator is stable; that is, it has strictly negative eigenvalues.  

We are currently studying the effect of different weighting functions on the operator spectrum for a future theory paper.  We are also looking into how other choices, like neighbor selection and centering of the interpolant, modify the spectrum.  Our current method produces a Laplacian operator with eigenvalues that depend on the inverse of the smallest volume fractions.  The operator spectrum also contains a small cluster of eigenvalues with non-negligible imaginary components.  We hope to alleviate both problems in the future.    Finally, we plan to study the effect of using standard regular stencils in the interior for the theory paper.

The method described in this paper is applicable to other problems in div-flux form.  For example, we are working on a fourth-order method for advection-diffusion.  For this problem, we are using an upwind weighting system.  Another application for the future is the variable coefficient Poisson's equations for smoothly varying coefficients.

%\input{appendix}

%\ack
% The author wishes to thank

%\clearpage

\bibliography{poisson}

\begin{thebibliography}{10}

\bibitem{Adams-03a}
M.~F. Adams.
\newblock Algebraic multigrid methods for constrained linear systems with
  applications to contact problems in solid mechanics.
\newblock {\em Numerical Linear Algebra with Applications}, 11(2-3):141--153,
  2004.

\bibitem{AftosmisBergerMelton}
M.~J. Aftosmis, M.~J. Berger, and J.~S. Saltzman.
\newblock Robust and efficient cartesian mesh generation for component-base
  geometry.
\newblock {\em AIAA Journal}, 36(6):952--960, June 1998.

\bibitem{petsc-user-ref}
S.~Balay, S.~Abhyankar, M.~F. Adams, J.~Brown, P.~Brune, K.~Buschelman,
  V.~Eijkhout, W.~D. Gropp, D.~Kaushik, M.~G. Knepley, L.~C. McInnes, K.~Rupp,
  B.~F. Smith, and H.~Zhang.
\newblock {PETS}c users manual.
\newblock Technical Report ANL-95/11 - Revision 3.5, Argonne National
  Laboratory, 2014.

\bibitem{petsc-web-page}
S.~Balay, S.~Abhyankar, M.~F. Adams, J.~Brown, P.~Brune, K.~Buschelman,
  V.~Eijkhout, W.~D. Gropp, D.~Kaushik, M.~G. Knepley, L.~C. McInnes, K.~Rupp,
  B.~F. Smith, and H.~Zhang.
\newblock {PETS}c {W}eb page.
\newblock http://www.mcs.anl.gov/petsc, 2014.

\bibitem{petsc-efficient}
S.~Balay, W.~D. Gropp, L.~C. McInnes, and B.~F. Smith.
\newblock Efficient management of parallelism in object oriented numerical
  software libraries.
\newblock In E.~Arge, A.~M. Bruaset, and H.~P. Langtangen, editors, {\em Modern
  Software Tools in Scientific Computing}, pages 163--202. Birkh{\"{a}}user
  Press, 1997.

\bibitem{BRENNERFEMBOOK}
S.~Brenner and R.~Scott.
\newblock {\em The Mathematical Theory of Finite Element Methods}.
\newblock Springer, New York, 2007.

\bibitem{slepc-users-manual}
C.~Campos, J.~E. Roman, E.~Romero, and A.~Tomas.
\newblock {SLEPc} users manual.
\newblock Technical Report DSIC-II/24/02 - Revision 3.3, D. Sistemes
  Inform\`atics i Computaci\'o, Universitat Polit\`ecnica de Val\`encia, 2012.

\bibitem{Cheng1999}
H.~Cheng, L.~Greengard, and V.~Rokhlin.
\newblock A fast adaptive multipole algorithm in three dimensions.
\newblock {\em JCP}, 155:468--496, 1999.

\bibitem{GibouFedkiw2005}
F.~Gibou and R.~Fedkiw.
\newblock A fourth order accurate discretization for the laplace and heat
  equations on arbitrary domains, with applications to the stefan problem.
\newblock {\em J. Comput. Phys.}, 202:577--601, 2005.

\bibitem{Greengard1996}
L.~Greengard and J.-Y. Lee.
\newblock A direct adaptive poisson solver of arbitrary order accuracy.
\newblock {\em J. Comput. Phys.}, 125:415--424, 1996.

\bibitem{Hernandez:2003:SSL}
V.~Hernandez, J.~E. Roman, and V.~Vidal.
\newblock {SLEPc}: {S}calable {L}ibrary for {E}igenvalue {P}roblem
  {C}omputations.
\newblock {\em Lecture Notes in Computer Science}, 2565:377--391, 2003.

\bibitem{Johansen:1998:EBPoisson}
H.~Johansen and P.~Colella.
\newblock A cartesian grid embedded boundary method for poisson's equation on
  irregular domains.
\newblock {\em J. Comput. Phys.}, 147(1):60--85, 1998.

\bibitem{johansenColella:1998}
H.~S. Johansen and P.~Colella.
\newblock A {Cartesian} grid embedded boundary method for {Poisson's} equation
  on irregular domains.
\newblock {\em J. Comput. Phys.}, 147(2):60--85, December 1998.

\bibitem{LevequeLing1994}
R.~LeVeque and Z.~Ling.
\newblock The immersed interface method for elliptic equations with
  discontinuous coefficients and singular sources.
\newblock {\em SIAM Journal of Numerical Analysis}, 31(4):1019--1044, 1994.

\bibitem{LEVEQUEBOOK}
R.~J. LeVeque.
\newblock {\em Numerical Methods for Conservation Laws}.
\newblock Birkhauser-Verlag, Basel, Boston, Berlin, 1990.

\bibitem{McKenneyFMM1994}
A.~McKenney, L.~Greengard, and A.~Mayo.
\newblock A fast poisson solver for complex geometries.
\newblock {\em J. Comput. Phys.}, 118:348--355, 1995.

\bibitem{PirzadehFEMGridGen2010}
S.~Pirzadeh.
\newblock Advanced unstructured grid generation for complex aerodynamic
  applications.
\newblock {\em AIAA Journal}, 48(5):904--915, 2010.

\bibitem{Romero:2014:PID}
E.~Romero and J.~E. Roman.
\newblock A parallel implementation of {Davidson} methods for large-scale
  eigenvalue problems in {SLEPc}.
\newblock {\em {ACM} Trans. Math. Software}, 40(2):13:1--13:29, 2014.

\bibitem{schwartzETAL:2006}
P.~Schwartz, M.~Barad, P.~Colella, and T.~Ligocki.
\newblock A {C}artesian grid embedded boundary method for the heat equation and
  {P}oisson's equation in three dimensions.
\newblock {\em Journal of Computational Physics}, 211(2):531--550, Jan. 2006.

\bibitem{EBGeometryPaper}
P.~Schwartz, J.~Percelay, T.~Ligocki, H.~Johansen, D.~T. Graves, D.~Devendran,
  P.~Colella, and E.~Ateljevich.
\newblock High accuracy embedded boundary grid generation using the divergence
  theorem.
\newblock Accepted by CAMCoS, 2014.

\end{thebibliography}
\bibliographystyle{abbrv}

\end{document}